\documentclass{gtart}
\def\ifplaintex{\expandafter\ifx\csname documentclass\endcsname\relax}

\ifplaintex 
\else
\fi

\expandafter\ifx\csname beginpicture\endcsname\relax
\expandafter\ifx\csname documentclass\endcsname\relax
\input pictex \else
\input prepictex \input pictex \input postpictex \fi\fi

\def\gt{{\mathsurround=0pt\it $\cal G\mskip-2mu$eometry \&\ 
$\cal T\!\!$opology}}        %  journal title in recommended style

\def\gtp{{\mathsurround=0pt\it $\cal G\mskip-2mu$eometry \&\ 
$\cal T\!\!$opology $\cal P\!$ublications}}  % GT publications

\def\lognumber#1{\def\thelognumber{#1}}
\def\volumenumber#1{\def\thevolumenumber{#1}}
\def\papernumber#1{\def\thepapernumber{#1}}
\def\volumeyear#1{\def\thevolumeyear{#1}}

\def\pagenumbers#1#2{\def\startpage{#1}\def\finishpage{#2}}
\def\published#1{\def\publishdate{#1}}
\def\proposed#1{\def\theproposer{#1}}
\def\seconded#1{\def\theseconders{#1}}
\def\received#1{\def\receiveddate{#1}}

\def\accepted#1{\def\accepteddate{#1}}

\def\coverauthors#1{\def\thecoverauthors{#1}}
\def\asciiauthors#1{\def\theasciiauthors{#1}}
\def\asciiaddress#1{\def\theasciiaddress{#1}}
\def\asciiemail#1{\def\theasciiemail{#1}}
\long\def\asciiabstract#1{\long\def\theasciiabstract{#1}}
\def\asciikeywords#1{\def\theasciikeywords{#1}}
\def\shortauthors#1{\def\theshortauthors{#1}}

\let\\\par\let\thelognumber\relax
\let\thevolumenumber\relax\let\thepapernumber\relax
\let\thevolumeyear\relax\let\thesamplenumber\relax\let\startpage\relax
\let\finishpage\relax\let\publishdate\relax\let\receiveddate\relax
\let\reviseddate\relax\let\accepteddate\relax\let\theasciititle\relax
\let\theasciiauthors\relax\let\theasciiaddress\relax
\let\theasciiabstract\relax\let\theasciikeywords\relax
\let\theasciiemail\relax\let\theshortauthors\relax\let\theshorttitle\relax
\let\thecoverauthors\relax

\long\def\maketitlep{   % start of definition of \maketitlep

\count0=\startpage

\gt\hfill      %   Journal title (top left) 
\beginpicture
\setcoordinatesystem units <0.33truein, 0.33truein> point at 2.2 0.9
\setplotsymbol ({$\cal G$})
\plotsymbolspacing=9truept
\circulararc 315 degrees from 0 1 center at 0 0
\setplotsymbol ({$\cal T$})
\circulararc 315 degrees from 1 -1 center at 1 0
\endpicture
\break
{\small\ifx\thesamplenumber\relax % sample?  
Volume \else Sample
\fi\thevolumenumber\ (\thevolumeyear)
\startpage--\finishpage\nl
Published: \publishdate}
\vglue 0.5truein plus 0.4fil minus 0.1truein

{\parskip=0pt\leftskip 0pt plus 1fil\def\\{\par\smallskip}{\ifplaintex\large
\else\Large\fi\bf\thetitle}\par\medskip}   

\vglue 0pt plus 0.1fil 

{\parskip=0pt\leftskip 0pt plus 1fil\def\\{\par}{\sc\theauthors}
\par\medskip}

\vglue 0pt plus 0.1fil 

{\small\parskip=0pt\let\newline\\
{\leftskip 0pt plus 1fil\def\\{\par}{\sl\theaddress}\par}
\expandafter\ifx\theemail\relax    % email address?
\relax\else\vglue 5pt plus 0.02fil minus 2pt\def\\{\stdspace{\rm 
and}\stdspace} 
\cl{Email:\stdspace\tt\theemail}\fi
\ifx\theurl\relax                  % URL given?
\relax\else\vglue 5pt plus 0.02fil minus 2pt\def\\{\stdspace{\rm 
and}\stdspace}
\cl{URL:\stdspace\tt\theurl}\fi\par}

\vglue 7pt plus 0.3fil minus 3pt

{\bf Abstract}
\vglue 5pt plus 0.1fil minus 2pt

\theabstract

\vglue 7pt plus 0.3fil minus 3pt

{\bf AMS Classification numbers}\quad Primary:\quad \theprimaryclass

Secondary:\quad \thesecondaryclass

\vglue 5pt plus 0.3fil minus 2pt

{\bf Keywords}\quad \thekeywords

\vglue 10pt plus 0.5fil minus 5pt

{\small  Proposed: \theproposer\hfill Received: \receiveddate\nl
Seconded: \theseconders\hfill 
\ifx\reviseddate\relax                         % paper revised?
Accepted: \accepteddate                        % no
\else
Revised: \reviseddate                          % yes
\fi}
\eject
}       %  end of definition of \maketitlep

\let\maketitlepage\maketitlep
\let\maketitle\maketitlepage

\font\phead=cmsl9 scaled 950
\font=cmsl9 scaled 1050
\font\pnum=cmbx10 scaled 913
\font\lnum=cmbx10 
\font\pfoot=cmsl9 scaled 950
\font=cmsl9 scaled 1050
\ifplaintex
\headline{\vbox to 0pt{\vskip -4.5mm\line{\small\phead\ifnum
\count0=\startpage ISSN 1364-0380 (on line)
1465-3060 (printed) \hfill {\pnum\folio}\else\ifodd\count0\def\\{ }% 
\ifx\theshorttitle\relax\thetitle\else\theshorttitle\fi\hfill{\pnum\folio}
\else\def\\{ and }{\pnum\folio}\hfill\ifx\theshortauthors\relax\theauthors
\else\theshortauthors\fi\fi\fi}\vss}}
\footline{\vbox to 0pt{\vglue 0mm\line{\small\pfoot\ifnum\count0=\startpage
\copyright\ \gtp\hfill\else
\gt, Volume \thevolumenumber\ (\thevolumeyear)\hfill\fi}\vss
}}
\else
\makeatletter
\def\@oddhead{{\small\lhead\ifnum\count0=\startpage ISSN 1364-0380 (on line)
1465-3060 (printed) \hfill {\lnum\number\count0}\else\ifodd\count0
\def\\{ }\ifx\theshorttitle\relax \thetitle \else\theshorttitle\fi\hfill
{\lnum\number\count0}\else\def\\{ and }{\lnum\number\count0}
\hfill\ifx\theshortauthors\relax 
\theauthors\else\theshortauthors\fi\fi\fi}}\def\@evenhead{\@oddhead}
\def\@oddfoot{\small\lfoot\ifnum\count0=\startpage\copyright\ \gtp\hfill\else
\gt, Volume \thevolumenumber\ (\thevolumeyear)\hfill\fi}
\def\@evenfoot{\@oddfoot}
\makeatother
\fi

\newwrite\gtoutfile
\long\gdef\makeheadfile{  %%% start of definition of \makeheadfile
{\def\\{, }\def\s{ }
\immediate\openout\gtoutfile head.xxx
\immediate\write\gtoutfile{To: math@arxiv.org}
\immediate\write\gtoutfile{Subject: put or rep NNNNN:pppp}
\immediate\write\gtoutfile{--text follows this line--}
\immediate\write\gtoutfile{Proxy-for: \ifx\theasciiauthors\relax
\theauthors\else\theasciiauthors\fi\s<\ifx\theasciiemail\relax\theemail\else\theasciiemail\fi>}
\immediate\write\gtoutfile{\noexpand\\}
\immediate\write\gtoutfile{Authors: \ifx\theasciiauthors\relax
\theauthors\else\theasciiauthors\fi}
{\def\\{ }\immediate\write\gtoutfile{Title: \ifx\theasciititle\relax
\thetitle\else\theasciititle\fi}}
\immediate\write\gtoutfile{Subj-class: GT or SG or MG etc}
\immediate\write\gtoutfile{MSC-class: \theprimaryclass\ifx\thesecondaryclass\relax\else, \thesecondaryclass\fi}
\immediate\write\gtoutfile{Journal-ref: Geom. Topol. \thevolumenumber
(\thevolumeyear) \startpage-\finishpage}
\immediate\write\gtoutfile{Comments: Published by Geometry and Topology at}
\immediate\write\gtoutfile{\s\s http://www.maths.warwick.ac.uk/gt/GTVol\thevolumenumber/paper\thepapernumber.abs.html}
\immediate\write\gtoutfile{\noexpand\\}
\immediate\write\gtoutfile{}
\ifx\theasciiabstract\relax
\immediate\write\gtoutfile{\theabstract}\else
\immediate\write\gtoutfile{\theasciiabstract}\fi
\immediate\write\gtoutfile{}
\immediate\write\gtoutfile{\noexpand\\}
\immediate\write\gtoutfile{}
\immediate\closeout\gtoutfile}}  %%% end of definition of \makeheadfile

\def\maketitlepage{\maketitlep\makeheadfile}
\let\maketitle\maketitlepage

\lognumber{294}
\volumenumber{7}\papernumber{15}
\volumeyear{2003}
\pagenumbers{511}{568}
\received{6 December 2002}
\accepted{8 June 2003}
\published{18 August 2003}
\proposed{Shigeyuki Morita}
\seconded{Ralph Cohen, Vaughan Jones}

\usepackage{amsmath,amssymb,amscd,epsf}

\swapnumbers
\newtheorem{subsubtheorem}[subsubsection]{}
\def\subsubth#1{\subsubtheorem{\bf#1}\rm\qua\ignorespaces}
\newtheorem{subtheorem}[subsection]{}
\def\subth#1{\subtheorem{\bf#1}\rm\qua\ignorespaces}

\def\S{section }
\let\cal\mathcal

\def\del{\partial}

\def\a{\alpha}
\def\b{\beta}

\def\l{\lambda}
\def\s{\sigma}

\def\th{\theta}
\def\Cal{\cal}

\def\ZZ{{\mathbb Z}/2{\mathbb Z}}

\def\nn{\nonumber}

\def\Zz{{\mathbb Z}/2{\mathbb Z}}

\def\R{\mathbb R}
\def\Z{\mathbb Z}
\def\SS{{\mathbb S}}
\def\Sn{{\mathbb S}_{n}}
\def\Snn{{\mathbb S}_{n+1}}

\def\lr{\stackrel{\leftrightarrow}{d}}
\def\ldel{\stackrel{\leftarrow}{\del}}
\def\rdel{\stackrel{\rightarrow}{\del}}
\def\rr{\stackrel{\rightarrow}{d}}

\def\darc{\mathcal{DARC}}
\def\Loop{\mathcal{LOOP}}

\def\Arc{\mathcal{ARC}}
\def\vardel{\delta}
\def\del{\partial}
\def\nn{\nonumber}
\def\s{\sigma}
\def\text{\mbox}

\begin{document}
\title{Arc Operads and Arc Algebras}

\author{Ralph M Kaufmann\\Muriel Livernet\\R\thinspace C Penner}
\coverauthors{Ralph M Kaufmann\\Muriel Livernet\\R\noexpand\thinspace C Penner}
\shortauthors{Ralph M Kaufmann, Muriel Livernet and R\thinspace C Penner}
\asciiauthors{Ralph M Kaufmann, Muriel Livernet and RC Penner}

\address{{\rm RMK:}\qua Oklahoma State University, Stillwater, USA\\
{\rm and}\qua Max--Planck--Institut f\"ur Mathematik, Bonn, Germany\\
\smallskip\\{\rm ML:}\qua Universit\'e Paris 13, France\\
\smallskip\\{\rm RCP:}\qua University of Southern California, Los
Angeles, USA\\\smallskip\\{\rm Email:}\qua {\tt ralphk@mpim-bonn.mpg.de,
livernet@math.univ-paris13.fr\\rpenner@math.usc.edu}}

\asciiaddress{RMK: Oklahoma State University, Stillwater, USA\\
and Max-Planck-Institut fur Mathematik, Bonn, Germany\\
ML: Universite Paris 13, France\\
RCP: University of Southern California, Los
Angeles, USA}

\asciiemail{ralphk@mpim-bonn.mpg.de,
livernet@math.univ-paris13.fr,
rpenner@math.usc.edu}

\begin{abstract}
Several topological and homological operads
based on families of projectively weighted arcs in bounded surfaces are
introduced and studied.  The spaces underlying the basic operad are
identified with open subsets of a combinatorial compactification due to
Penner of a space closely related to Riemann's moduli space.  Algebras
over these operads are shown to be Batalin--Vilkovisky algebras, where the
entire BV structure is realized simplicially.  Furthermore, our basic
operad contains the cacti operad up to homotopy. New operad structures
on the circle are classified and combined with the basic operad to
produce geometrically natural extensions of the algebraic structure of BV
algebras, which are also computed.
\end{abstract}

\asciiabstract{Several topological and homological operads based on
families of projectively weighted arcs in bounded surfaces are
introduced and studied. The spaces underlying the basic operad are
identified with open subsets of a compactification due to Penner of a
space closely related to Riemann's moduli space. Algebras over these
operads are shown to be Batalin-Vilkovisky algebras, where the entire
BV structure is realized simplicially. Furthermore, our basic operad
contains the cacti operad up to homotopy, and it similarly acts on the
loop space of any topological space. New operad structures on the
circle are classified and combined with the basic operad to produce
geometrically natural extensions of the algebraic structure of BV
algebras, which are also computed.}

\asciikeywords{Moduli of Surfaces, Operads, Batalin-Vilkovisky algebras}



\maketitlepage

\section*{Introduction}
With the rise of string theory and its companions topological and
conformal field theories, there has been a very fruitful
interaction between physical ideas and topology. One of its
manifestations is the renewed interest and research activity in
operads and moduli spaces, which has uncovered many new
invariants, structures and many unexpected results.

The present treatise adds to this general understanding. The
starting point is Riemann's moduli space, which provides the basic
operad underlying conformal field theory. It has been extensively
studied and gives rise to Batalin--Vilkovisky structures and
gravity algebras. Usually one is also interested in
compactifications or partial compactifications of this space which
lead to other types of theories, as for instance, the
Deligne--Mumford compactification corresponds to cohomological
field theories as are present in Gromov--Witten theory and quantum
cohomology. There are other types of compactifications however
that are equally interesting like the new compactification of Losev and Manin
\cite{LM} and the combinatorial compactification of Penner \cite{P1}.

The latter compactification, which we will employ in this paper, 
is for instance well suited to handle matrix theory and
constructive field theory. The combinatorial nature of the
construction has the great benefit that all calculations can be
made explicit and often have very nice geometric
descriptions. The combinatorial structure of this space is
described by complexes of embedded {weighted} arcs on surfaces, 
upon which we shall comment further below.
These complexes themselves are of extreme interest as they provide
a very general description of surfaces with extra structure which
can be seen as acting via cobordisms. The surfaces devoid of arcs
give the usual cobordisms defining topological field theories by
Atiyah and Segal. Taking the arc complexes as the basic objects, we
construct several operads and suboperads out of them, by
specializing the general structure to  more restricted ones,
thereby controlling the ``information carried by the cobordisms''.

The most interesting one of these is the cyclic operad $\Arc$ that
corresponds to surfaces with weighted arcs such that all boundary
components have arcs emanating from them. We will call algebras
over the corresponding homology operad Arc algebras. The imposed
conditions correspond to open subsets of the combinatorial
compactification which contain the operators that yield the
Batalin--Vilkovisky structure. Moreover, this structure is very natural
and is topologically explicit in a clean and beautiful way.
 This operad in genus zero with
no punctures can also be embedded into an operad in arbitrary genus with
an arbitrary number of punctures with the help of a
genus generator and a puncture operator, which ultimately should
provide a good control for stabilization.

\begin{figure}[ht!]
\epsfxsize = \textwidth
\epsfbox{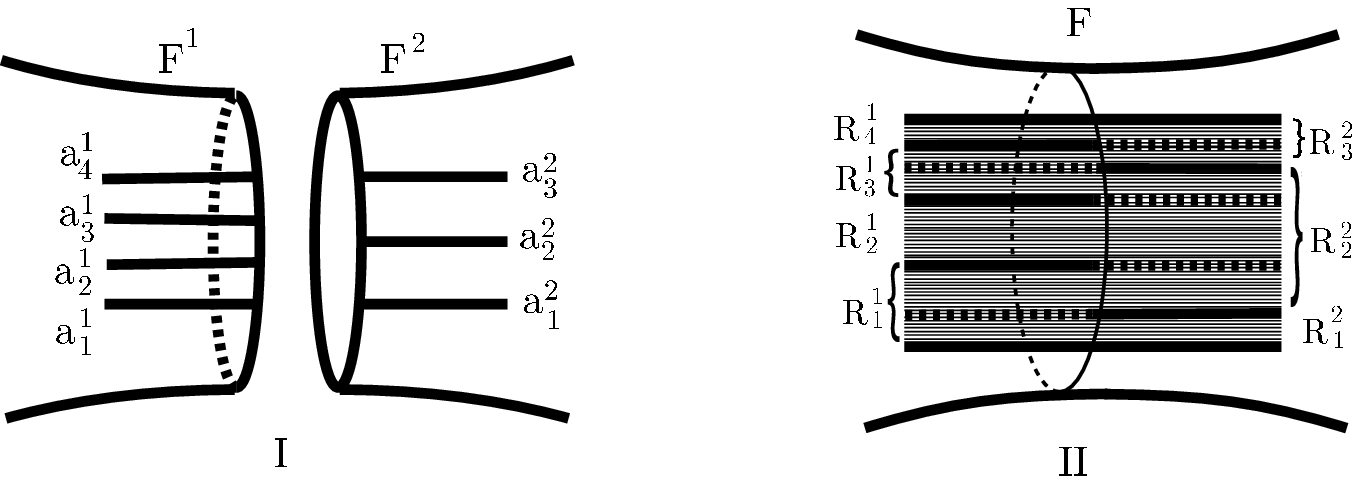}

\caption{\label{introfig}
Gluing bands from weighted arcs:}
\cl{\small I, arc families in two surfaces,\qua II, combining bands.}
\end{figure}

A number of operads of current interest \cite{CS,J,V,MSS} are governed by a composition which, in effect, depends
upon combining families of arcs in surfaces.  Our operad composition on $\Arc$ depends upon
an explicit method of
combining families of {\sl weighted} arcs in surfaces, that is, each component arc of the family is
assigned a positive real number; we
next briefly describe this composition.  Suppose that
$F^1,F^2$ are surfaces with distinguished boundary components $\partial ^1\subseteq F^1, \partial
^2\subseteq F^2$.  Suppose further that each surface comes equipped with a properly embedded family of
arcs, and let
$a_1^i,\ldots ,a_{p^i}^i$ denote the arcs in
$F^i$ which are incident on $\partial ^i$, for $i=1,2$, as illustrated in part I of figure \ref{introfig}.
Identify $\partial ^1$ with $\partial ^2$ to produce a surface $F$.  We wish to furthermore combine
the arc families in $F^1,F^2$ to
produce a corresponding arc family in $F$, and there is evidently no well-defined way to achieve this
without making further choices or imposing further conditions on the arc families (such as $p^1=p^2$).

Our additional data required for gluing is given by an assignment of one real number, a weight, to each arc
in each of the arc families.  
The weight $w_j^i$ on $a_j^i$ is interpreted geometrically as the height of a rectangular band
$R_j^i=[0,1]\times [-{{w_j^i}/ 2}, {{w_j^i}/ 2}]$ whose core $[0,1]\times\{ 0\}$ is identified with $a_j^i$, for $i=1,2$ and
$j=1,\ldots ,p^i$. We shall assume that
$\sum _{j=1}^{p^1} w^1_j=\sum _{j=1}^{p^2} w^2_j$ for simplicity, so that the total height of all the bands
incident on $\partial ^1$ agrees with that of $\partial ^2$; in light of this assumption, the bands in
$F^1$ can be sensibly attached along $\partial ^1$ to the bands in $F^2$ along $\partial ^2$ to
produce a collection of bands in the surface $F$ as illustrated in part II of figure \ref{introfig};
notice that the horizontal edges of the rectangles $\{ R_j^1\} _1^{p^1}$ decompose the rectangles 
$\{ R_j^2\} _1^{p^2}$ into sub-rectangles and conversely.  The resulting family of sub-bands, in turn,
determines a weighted family of arcs in $F$, one arc for each sub-band with a weight given by the width
of the sub-band; thus, the weighted arc family in $F$ so produced depends upon the weights in a
non-trivial but combinatorially explicit way.  This describes the basis of our gluing operation on families of weighted arcs, which is
derived from the theory of train tracks and partial measured foliations \cite{PH}.  In fact, the simplifying
assumption that the total heights agree is obviated by considering not weighted families of arcs, but
rather weighted families of arcs modulo the natural overall homothetic action of
${\mathbb R}_{>0}$--so-called projectively weighted arc families; given projectively weighted arc families, we may
de-projectivize in order to arrange that the simplifying assumption is in force (assuming that $p^1\neq
0\neq p^2$), perform the construction just described, and finally re-projectivize.  We shall prove in
\S 1 that this construction induces a well-defined operad composition on suitable classes of projectively
weighted arcs.

The resulting operad structure provides a useful framework as is evidenced by the fact that there is an embedding of the
cactus operad of Voronov \cite{V} into this operad as a suboperad, as we
discovered. This shows how to view the Chas--Sullivan \cite{CS} 
string topology, which was the inspiration for our
analysis and to which \S2 owes  obvious intellectual gratitude,
inside the combinatorial model. In particular, we recover
in this way a surface description of the cacti which in turn can
be viewed as a reduction of the surface structure in a very
precise way. In fact, there is a reduction of any surface with
arcs to a configuration in the plane, which is not necessarily a
cactus and may have a much more complex structure. For the
surfaces whose reduction is a cactus in the sense of \cite{V}, one
retains the natural action on loop spaces by forgetting some
of the internal topological structure of the surface but
keeping an essential part of the information carried by the
arcs. For more general configurations, a similar but more
complicated structure is expected.

One virtue of viewing the cacti as a surface instead of as a
singular level set is that in this way the branching behavior is
nicely depicted while the ``singularities of one--dimensional
Feynman graphs do not appear'' as Witten has pointed out many
times.  Furthermore, there are natural suboperads
governing the Gerstenhaber and the BV structure, where the latter
suboperad is generated by the former and the 1-ary operation of
$\Arc$. On the level of homology, these operation just add one class,
that of the BV operator. For cacti, the situation is much more
complicated and is given by a bi--crossed product \cite{K}.

It is conjectured \cite{P1} that the full arc complex of a surface with boundary
is spherical and
thus will not carry much operadic information on the homological level. 
The
sphericity of the full compactification can again be compared to
the Deligne--Mumford situation where the compactification in the
genus zero case leaves no odd cohomology and in a Koszul dual way
is complementary to the gravity structure of the open moduli
space. In general the idea of obtaining suboperads by
imposing certain conditions can be seen as parallel to the
philosophy of Goncharov--Manin \cite{GM}. 
The Sphericity Conjecture \cite{P1} would 
provides the basis for a calculation of the homology of the
operads under consideration here, a task that will be undertaken in
\cite{KLP}.

One more virtue of having concrete surfaces is that we can also
handle additional structures on the objects of the corresponding
cobordism category, viz.\ the boundaries. We incorporate these
ideas in the form of direct and semi--direct products of our
operads with operads based on circles. These give geometrically natural
extensions of the algebraic structure of BV algebras.  This view is also inherent
in \cite{K} where the cacti operads are decomposed into bi--crossed
products. The operads built on circles which we consider have the
geometric meaning of marking additional point in the boundary of a surface. Their algebraic
construction, however, is not linked to their particular
presentation but rather relies on the algebraic structures of
$S^1$ as a monoid. They are therefore also of
independent algebraic interest.

The paper is organized as follows: In \S1 we review the salient
features of the combinatorial compactification of moduli spaces and the 
Sphericity Conjecture. Using this background, we define the operadic
products on surfaces with weighted arcs that underlie all of our
constructions. In \S2 we uncover the
Gerstenhaber and Batalin--Vilkovisky  structures of our operad on
the level of chains with explicit chain homotopies.  These chain homotopies
manifestly show the symmetry of these equations. The BV--operator
is given by the unique, up to homotopy, cycle of arc families on
the cylinder. Section 3 is devoted to the inclusion of the cacti operad
into the arc operad both in its original version \cite{CS,V} as well as
in its spineless version \cite{K}. This explicit map, called framing,
also shows that the spineless cacti govern the Gerstenhaber
structure while the Voronov cacti yield BV; this fact can also
be read off from the explicit calculations of \S2 for the respective
suboperads of $\Arc$, which are also defined in this section. The
difference between the two sub--operads are operations
corresponding to a Fenchel--Nielsen type deformation of the
cylinder, i.e., the 1-ary operation of $\Arc$. 
A  general construction which forgets the topological structure of
the underlying surface and retains a collection of parameterized
loops in the plane with incidence/tangency conditions dictated by
the arc families is contained in \S4. Using this partial forgetful
 map we  define the action of a suboperad of
$\Arc$ on loop spaces of manifolds. 
For the particular
suboperads corresponding to cacti, this operation is inverse to the
framing of \S3. In \S5 we define several direct and semi--direct
products of our operad with cyclic and non--cyclic operads built
on circles. One of these products yields for instance dGBV algebras.
The operads built on circles which we utilize are provided by our
analysis of this type of operad in the Appendix. The approach is a
classification of all operations that are linear and local in the
coordinates of the components to be glued. The results are not
contingent on the particular choice of $S^1$, but only on certain
algebraic properties like being a monoid and are thus of a more
general nature.

\rk{Acknowledgments}
The first author wishes to thank
 the IH\'ES in Bures--sur--Yvette
and the Max--Planck--Institut
f\"ur Mathematik in Bonn, where some of the work was carried out,
for their hospitality. The second author would like to thank USC 
for its hospitality.
The financial support of the NSF under grant DMS\#0070681 is also
acknowledged by the first author as is the support of NATO by
the second author.

\section{The operad of weighted arc families in surfaces}

It is the purpose of this section to define our basic topological operad.  The idea of our operad
composition on projectively weighted arc families was already described in the Introduction.  The technical details 
of this are somewhat involved, however, and we next briefly survey the material in this section.
To begin, we define weighted arc families in surfaces together with several different geometric models of the common 
underlying combinatorial structure.  The collection of projectively weighted arc families in a fixed surface is found
to admit the natural structure of a simplicial complex, which descends to a CW decomposition of the quotient of this simplicial complex
by the action of the  pure mapping class group of the surface.  This CW complex had arisen in Penner's
earlier work, and we discuss its relationship to Riemann's moduli space in an extended remark.  The
discussion in the Introduction of bands whose heights are given by weights on an arc family is formalized by
Thurston's theory of partial measured foliations, and we next briefly recall the salient details of this
theory.  The spaces underlying our topological operads are then introduced; in effect, we consider arc
families so that there is at least one arc in the family which
is incident on any given boundary component.  In this setting, we define a collection of abstract measure
spaces, one such space for each boundary component of the surface, associated to each appropriate
projectively weighted arc family.  We next define the operad composition of these projectively weighted arc
families which was sketched in the Introduction and prove that this composition is well-defined.  Our basic
operads are then defined, and several extensions of the construction are finally discussed.

\subsection{Weighted arc families}

\subsubsection{Definitions}
\label{defn1}

Let $F=F_{g,r}^{s}$ be a fixed oriented topological surface
of genus $g\geq 0$ with $s\geq 0$ punctures and $r\geq 1$ boundary
components, where $6g-7+4r+2s\geq
0$, so the boundary $\partial F$ of $F$ is necessarily non-empty by
hypothesis.  Fix an enumeration
$\partial _1,\partial _2,
\ldots ,
\partial _r$ of the boundary components of $F$ once and for all.
In each boundary
component $\partial _i$ of $F$, choose once and for all a closed
arc
$W_i\subset\partial _i$, called a {\it window}, for each $i=1,2,\ldots , r$.

The {\it pure mapping class group} $PMC=PMC(F)$ is the group of
isotopy
classes of all
orientation-preserving homeomorphisms of $F$ which fix each 
$\partial _i-W_i$ pointwise 
(and fix each $W_i$ setwise), for each $i=1,2,\ldots , r$.

Define an {\it essential arc} in $F$ to be an embedded path $a$
in $F$ whose endpoints lie among the windows, where we demand that $a$
is not isotopic rel endpoints to a path lying in $\partial F$.  Two
arcs
are said to be {\it parallel} if there is an isotopy between them
which fixes each $\partial _i -W_i$ pointwise (and fixes each $W_i$
setwise) for
$i=1,2,\ldots ,r$. 

An {\it arc family} in
$F$ is the isotopy class of an unordered collection of disjointly embedded
essential arcs in
$F$, no two of which are parallel.  Thus, there is a well-defined action of $PMC$ on arc families.

\subsubsection{Several models for arcs}
We shall see that collections of arc families in a fixed surface
lead to a natural simplicial complex, an open subset of which will form the topological spaces
underlying our basic operads.  Before turning to this discussion, though, let us briefly analyze
the role of the windows in the definitions above and describe several geometric models of the
common underlying combinatorics of arc families in bounded surfaces in order to put this role
into perspective.

\begin{figure}[ht!]
\epsfxsize = \textwidth
\epsfbox{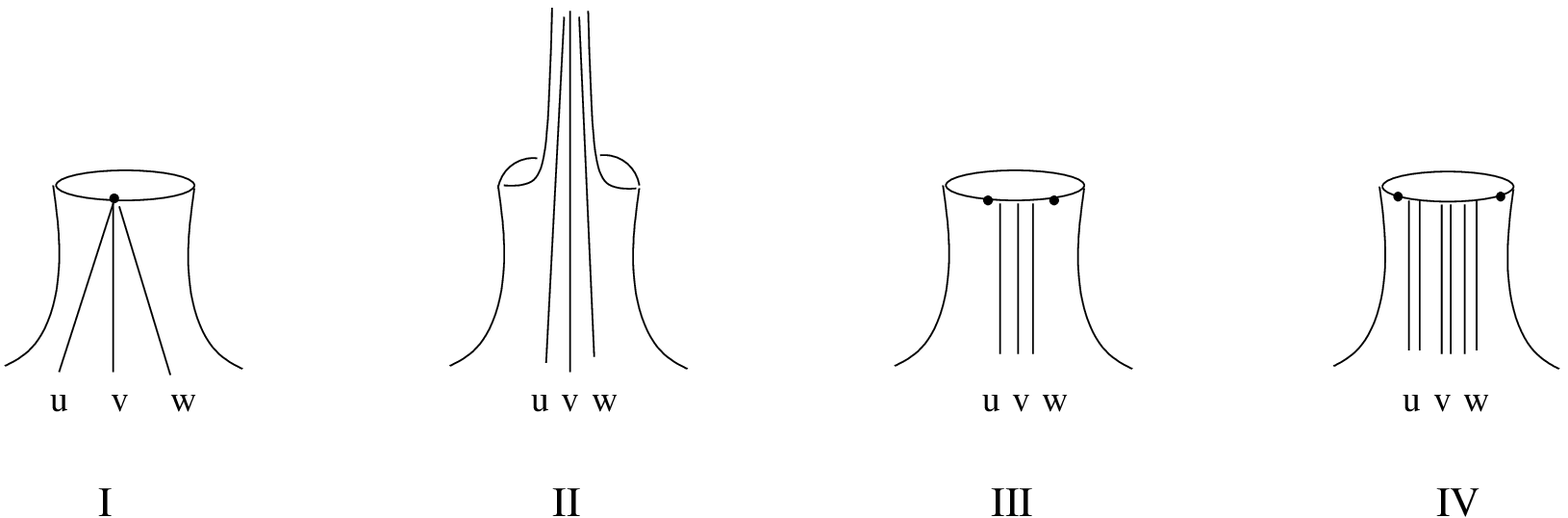}
\caption{\label{boundary}
I, arcs running to a point on the boundary,\qua II, arcs running
to a point at infinity,\qua
III, arcs in a window,\qua IV, bands in a window}
\end{figure}

For the first such model, let us choose a distinguished point $d_i\in\partial _i$, for $i=1,2,\ldots
,r$, and consider the space of all complete finite-area metrics on $F$ of constant Gauss curvature
$-1$ (so-called ``hyperbolic metrics'') so that each $\partial ^\times_i=\partial _i-\{ d_i\}$
is totally geodesic (so-called ``quasi hyperbolic metrics'') on $F$.  To explain this, consider a
hyperbolic metric with geodesic boundary on a once-punctured annulus
$A$ and the simple geodesic arc $a$ in it asymptotic in both directions to the puncture; the induced
metric on a component of $A-a$ gives a model for the quasi hyperbolic structure on $F^\times =F-\{
d_i\} _1^r$ near $\partial _i^\times$.  The first geometric model for an arc family $\alpha$ in $F$ is
a set of disjointly embedded geodesics in $F^\times$, each component of which is asymptotic in both directions to
some distinguished point $d_i$; see part II of figure \ref{boundary}.  

In the homotopy class of each $\partial _i$, there is a unique geodesic $\partial _i^*\subset
F^\times$.  Excising from $F-\cup\{ \partial_i^*\}_1^r$ 
any component which contains a point of
$\partial ^\times$, we obtain a hyperbolic structure on the surface $F^*\subseteq F^\times$ with
geodesic boundary (where in the special case of an annulus, $F^*$ collapses to a circle). 
Taking $\alpha\cap F^*$, we find a collection of geodesic arcs connecting boundary components
(where in the special case of the annulus, we find two points in the circle).  

This is our second geometric model for arc families.  We may furthermore choose a distinguished
point $p_i\in\partial _i^*$ and a regular neighborhood $U_i$ of $p_i$ in $\partial _i^*$, for
$i=1,2,\ldots ,r$.  Provided
$p_i\notin\alpha$, we may take $U_i$ sufficiently small that $U_i\cap
\alpha=\emptyset$, so the arc $V_i=\partial _i^*-U_i$ forms a natural ``window'' containing
$\alpha\cap\partial _i^*$.  There is then an ambient isotopy of $F^*$ which shrinks each window $V_i$
down to a small arc $W_i\subseteq\partial _i^*$, under which $\alpha$ is transported to a family of
(non-geodesic) arcs with endpoints in the windows $W_i$.  In case $p_i$ does lie in $\alpha$, then
let us simply move $p_i$ a small amount in the direction of the natural orientation 
(as a boundary component of $F^*$) along $\partial _i^*$ and perform the same construction; see part III of
figure \ref{boundary}.

This leads to our final geometric model of arc families, namely, the model defined in \ref{defn1}
for the purposes of this paper, of arcs in a bounded surface
with endpoints in windows.  This third model is in the spirit of
train tracks and measured foliations (cf. \cite{PH}) as we shall see and is most convenient for describing the
operadic structure. 

\subsubsection{The arc complex}
Let us inductively build a simplicial complex $Arc'=Arc' (F)$ as
follows.  For each singleton arc family in $F$, there is a distinct vertex of $Arc'(F)$.  
Having
thus inductively constructed the $(k-1)$-skeleton of $Arc'$ for $k\geq 1$, let us add a
$k$-simplex
$\sigma (\alpha ')$ to $Arc'$ for each arc family $\alpha '$ in $F$ of cardinality $k+1$. 
The simplicial structure on
$\sigma (\alpha ')$ itself is the natural one, where faces correspond to sub-arc families of
$\alpha '$, and we may therefore identify the proper faces of $\sigma (\alpha ')$ with
simplices in the
$(k-1)$-skeleton of $Arc'$.  Adjoining
$k$-simplices in this manner for each such arc family $\alpha '$ of cardinality $k+1$ defines
the
$k$-skeleton of $Arc'$.  This completes the inductive definition of the simplicial complex
$Arc'$.

$Arc'$ is thus
a
simplicial complex, upon which
$PMC=PMC(F)$ acts continuously, and we define the quotient topological
space to be
$Arc=Arc(F)=Arc'(F)/PMC(F)$.  If $\alpha '=\{a_0,a_1,\ldots
,a_k\}\in Arc'$, then the arcs $a_i$ come in a canonical linear
ordering.  Namely, the orientation of $F$ induces an orientation on
each window, and traversing the windows in the order
$W_1,W_2,\ldots ,W_r$ in these orientations, one first encounters an
endpoint of the arcs
$a_0,a_1,\ldots , a_k$ in some order, which prescribes the claimed
linear
ordering.  Thus, cells in $Arc'$ cannot
have finite
isotropy in $PMC$, and the simplicial decomposition of $Arc'$ descends to a CW decomposition
of $Arc$ itself.

\subsubth{Sphericity Conjecture} \cite{P1}\qua
\sl Fix any surface $F=F_{g,r}^s$
\label{sphere}
with $6g-7+4r+2s\geq 0$.
Then $Arc(F_{g,r}^s)$ is piecewise-linearly
 homeomorphic to a sphere of dimension $6g-7+4r+2s$.
\rm

\medskip Recent work seems to provide a proof of this conjecture at least
in the case $g=0$ of planar surfaces as will be taken
up elsewhere.  It is worth emphasizing that the current paper is
independent of the Sphericity Conjecture, but this
gives a useful perspective on the relationship between Riemann's moduli
space and the operads studied here.

\subsubth{Remark}  \label{modulispace}
It is the purpose of this remark to explain the geometry
underlying
$Arc(F)$ which was uncovered in \cite{P1,P2}.  Recall the distinguished point
$p_i\in\partial _i$
chosen in our second geometric model for arc families.
The ``moduli space'' $M=M(F)$ of the surface $F$ with boundary is the
collection of all complete
finite-area metrics of constant Gauss curvature -1 with geodesic boundary,
together with a distinguished point
$p_i$ in each boundary component, modulo push forward by diffeomorphisms.
There is a natural action of ${\bf R}_+$
on $M$ by simultaneously scaling each of the hyperbolic lengths of the
geodesic boundary components, and we let
$M/{\bf R}_+$ denote the quotient.
The main result of \cite{P2} is that $M/{\bf R}_+$ is proper homotopy
equivalent to the complement of a codimension-two subcomplex $Arc_\infty
(F)$  of $Arc(F)$, where $Arc_\infty (F)$ corresponds
to arc families $\alpha$ so that some component of $F-\cup\alpha$ is other
than a polygon or once-punctured polygon.
It is remarkable to suggest that by adding to a space homotopy equivalent
to $M/{\bf R}_+$ a suitable simplicial complex we obtain a
sphere.  There is much known \cite{P1} about the geometry and combinatorics of
$Arc(F)$ and $Arc_\infty (F)$.
Furthermore, the Sphericity Conjecture should be useful in
calculating the homological operads of the arc operads.

\subsubsection{Notation} We shall always adopt the notation that if
$\alpha$ is the $PMC$-orbit of an arc family in $F$, then $\alpha '\in
Arc'$ denotes some chosen arc family representing $\alpha$.
Furthermore, it will sometimes be convenient to specify the components
$a_0,a_1,\ldots , a_k$, for $k\geq 0$, of an arc family
$\alpha '$ representing $\alpha\in Arc$, and we shall write simply
$\alpha '=\{ a_0,a_1,\ldots ,a_k\}$ in this case.  Since the linear
ordering on components of $\alpha '\in Arc'$ is invariant under the
action of $PMC$, it descends to a well-defined linear ordering on the
components underlying the $PMC$-orbit $\alpha$.

Of course, a $k$-dimensional cell in the CW decomposition of
$Arc$ is determined by the $PMC$-orbit of some arc family $\alpha '=\{
a_0,a_1,\ldots, a_k\}\in Arc'$, and (again relying on the
$PMC$-invariant linear ordering discussed previously) a point in the
interior of this cell is determined by the projective class of a
corresponding $(k+1)$-tuple
$(w_0,w_1,\ldots ,w_k)$ of positive real numbers.
As usual, we shall let
$[w_0:w_1:\cdots :w_k]$ denote affine coordinates on the projective
classes of corresponding non-negative real $(k+1)$-tuples, and if
$w_i=0$,
for
$i$ in a proper subset
$I\subseteq\{0,1,\cdots , k\}$, then the point of $Arc$ corresponding
to $(\alpha ',[w_0:w_1:\cdots :w_k])$ is identified with the point in
the cell corresponding to $\{a_j\in\alpha ':j\notin I\}$ with
projective
tuple gotten by deleting all zero entries  of
$[w_0:w_1:\cdots :w_k]$.

To streamline the notation, if 
$\alpha '=\{ a_0,a_1,\ldots ,a_k\}$ is in
$Arc'$ and if $(w_0,w_1,\ldots ,$ $w_k)\in {\mathbb R}_+^{k+1}$ is an
assignment of numbers, called {\it weights}, one weight to each
component
of
$\alpha '$, then
we shall sometimes suppress the weights by letting
$(\alpha ')$ denote the arc family $\alpha '$ with corresponding
weights
$w(\alpha ')= (w_0,w_1,\ldots ,w_k)$.  In the same manner, the
projectivized weights
on the components of $\alpha '$ will sometimes be suppressed, and
we let $[\alpha ']$ denote the arc family $\alpha '$ with projective
weights $w[\alpha ']=[w_0:w_1:\ldots :w_k]$.  These same notations
will be used for $PMC$-classes of arc families as well, so, for
instance,
$(\alpha )$ denotes the $PMC$-orbit $\alpha$ with weights $w(\alpha
)\in{\mathbb R}_+^{k+1}$, and a point in $Arc$ may be denoted
simply $[\alpha ]\in Arc(F)$, where the projective weight is given by
$w[\alpha ]$.

\subsection{Partial measured foliations}
Another point of view on elements of $Arc$, which is
useful in the subsequent constructions, is derived from the theory
of ``train tracks'' (cf. \cite{PH}) as follows.  If
$\alpha '=\{ a_0,a_1,\ldots ,a_k\}\in Arc'$
is given weights $(w_0,w_1,\ldots
,w_k)\in{\mathbb R}_+^{k+1}$, then we may regard $w_i$ as a transverse
measure
on $a_i$, for each $i=0,1,\ldots , k$ to determine a
``measured train track with
stops'' and corresponding ``partial measured foliation'', as considered in \cite{PH}.

For the convenience of the reader, we next briefly recall the salient
and elementary features of this construction.  Choose for the
purposes of
this discussion any complete Riemannian metric $\rho$ of finite area
on
$F$, suppose that each $a_i$ is smooth for $\rho$,
and consider for each $a_i$ the ``band'' $B_i$ in $F$ consisting of
all
points within $\rho$-distance $w_i$ of $a_i$.  If it is necessary,
scale
the metric $\rho$ to $\lambda\rho$, for $\lambda >1$, to guarantee
that
these bands are embedded, pairwise disjoint in $F$, and have their 
endpoints lying among the windows.  The band
$B_i$ about
$a_i$ comes equipped with a foliation by the arcs parallel to $a_i$
which
are a fixed $\rho$-distance to
$a_i$, and this foliation comes equipped with a transverse measure
inherited from $\rho$; thus,
$B_i$ is regarded as a rectangle of width
$w_i$ and some irrelevant length, for
$i=0,1,\ldots , k$.  The foliated and transversely measured bands
$B_i$,
for $i=0,1,\ldots , k$ combine to give a ``partial measured foliation'' of
$F$, that is, a foliation of a closed subset of $F$ supporting an invariant transverse
measure (cf. \cite{PH}).
The isotopy class in $F$ rel $\partial F$ of this
partial measured foliation is independent of the choice of metric
$\rho$.

Continuing to suppress the choice of metric $\rho$, for each
$i=1,2,\dots ,r$, consider
$\partial _i \cap \bigl (
\coprod _{j=0}^k B_j\bigr )$, which is empty if $\alpha$ does not meet $\partial _i$ and is
otherwise a collection of closed intervals
in $W_i$ with disjoint interiors.  Collapse to a point each component
complementary to the interiors of these intervals in $W_i$ to
obtain an interval, which we shall denote $\partial
_i(\alpha ')$. Each such interval $\partial _i(\alpha ')$ inherits an
absolutely continuous measure $\mu ^i$ from the transverse measures
on the
bands.

\subsubth{Definition} Given two representatives $(\alpha _1')$ and $(\alpha _2')$ of the
same\break
weight\-ed $PMC$-orbit $(\alpha )$, the respective measure spaces 
$(\partial _i(\alpha _1'), \mu ^i_1)$ and\break
$(\partial _i(\alpha _2'), \mu ^i_2)$ are canonically identified,
which
allows us to consider the measure space
$(\partial _i (\alpha ),
\mu ^i)$ of a weighted $PMC$-orbit $(\alpha )$ itself, which is called
the $i^{\rm th}$
{\it end\/} of $(\alpha )$, for each $i=1,2,\ldots ,r$.

\subsection{The spaces underlying the topological operad}

\subsubsection{Definitions}

An arc family $\alpha '$ in $F$ is said to be
{\it exhaustive} if
for each boundary component $\partial _i$, for $i=1,2,\ldots ,r$,
there is
at least one component arc in $\alpha '$ with its endpoints in the
window
$W_i$.  Likewise, a
$PMC$-orbit
$\alpha$ of arc families is said to be {\it exhaustive} if some
(that is, any) representing arc family $\alpha '$ is so.  
Define the topological spaces
\begin{eqnarray}
Arc_g^s(n)&=&\{ [\alpha ]\in
Arc(F_{g,n+1}^s): \alpha ~{\rm is~exhaustive}\} ,\nonumber\\
\widetilde{Arc}_g^s(n)&=&Arc_g^s(n)\times (S^1)^{n+1},\nonumber
\end{eqnarray}
which are the spaces that comprise our various families of topological operads.
In light of
Remark \ref{modulispace}, $Arc_g^s(n)$ is identified with an open
subspace of $Arc(F)$ properly containing a space homotopy equivalent to $M(F_{g,n+1}^s)/S^1$, and $Arc(F)$ is spherical in at least
the planar case.

Enumerating the boundary components of $F$ as $\partial _0,
\partial _1,\ldots ,\partial _n$ once and for all
(where $\partial _0$
will play a special role in the subsequent discussion) and letting
${\mathbb S}_p$
denote the $p^{\rm th}$ symmetric group, there is a natural
${\mathbb S}_{n+1}$-action on the labeling of boundary components which
restricts
to a natural
${\mathbb S}_n$-action on the boundary components labeled $\{ 1,2,\ldots ,
n\}$.
Thus, ${\mathbb S}_n$ and ${\mathbb S}_{n+1}$ act on $Arc_g^s(n)$, and extending by the
diagonal action of ${\mathbb S}_{n+1}$ on $(S^1)^{n+1}$, the symmetric groups
${\mathbb S}_n$
and ${\mathbb S}_{n+1}$ likewise act on $\widetilde{Arc}_g^s(n)$, where ${\mathbb S}_n$ by
definition acts trivially on the first coordinate in $(S^1)^{n+1}$.

Continuing with the definitions, if $[\alpha ]\times (t_0,t_1, \ldots ,
t_n)
\in\widetilde{Arc}_g^s(n)$, let us choose a corresponding
deprojectivization $(\alpha )$, fix some boundary component $\partial
_i$, for $i=0,1,\ldots ,n$, and let $m_i=\mu ^i(\partial _i(\alpha ))$
denote the total measure.  There is then a unique
orientation-preserving mapping
\begin{equation}
\label{circleboundary}
c_i^{(\alpha )}: \partial _i(\alpha )\to S^1
\end{equation}
which maps the (class of the) first
point of
$\partial _i(\alpha )$ in the orientation of the window $W_i$ to $0\in
S^1$, where the measure on the domain is  $\frac{2\pi}{m_i}~\mu ^i$
and on the range is the Haar measure on $S^1$. Of course, each
$c_i^{(\alpha )}$ is injective on the interior of $\partial _i(\alpha
)$ while $(c_i^{(\alpha )})^{-1} (\{ 0\} ) =\partial \bigl (\partial
_i(\alpha )\bigr )$.  Furthermore, if $(\alpha _1)$
and $(\alpha _2)$ are two different deprojectivizations of a common
projective class, then $c_i^{(\alpha _2)}\circ (c_i^{(\alpha
_1)})^{-1}$ extends continuously to the identity on $S^1$.

\subsubth{Remark}\label{remnu}
We do not wish to specialize to one or another particular case at this
stage, but if $\partial F$ comes equipped with an {\it a priori}
absolutely continuous measure (for instance, if
$\partial F$ comes equipped with a canonical
coordinatization, or if $F$ or $\partial F$ comes equipped with a
fixed
Riemannian metric), then we can identify $\partial _i (\alpha )$
with $\partial _i F$, for each $i=0,1,\ldots ,n$,
in the obvious manner (where, if the {\it a priori} total measure of
$\partial _i$ is
$M_i$, then we alter the Riemannian metric $\rho$
employed in the definition of the bands so that the $\mu _i$ total
measure
$m_i=\mu _i(\partial _i(\alpha ))$ agrees with $M_i$ for each
$i=0,1,\ldots ,n$ and finally collapse the endpoints of $\partial_i(\alpha)$ 
to a single point).

Thus, the idea is that the $i^{\rm th}$ end of $(\alpha )$ gives a
model
or coordinatization for $\partial _i$, but of course, altering the
underlying
$(\alpha )$ in turn alters the coordinatization $c_i^{(\alpha )}$ of
$\partial _i$.
The trivial product $\widetilde{Arc}_g^s(n)=Arc_g^s(n)\times
(S^1)^{n+1}$ leads in this way to a family of coordinatizations of
$\partial F$ which are twisted by $Arc_g^s(n)$.

\subsubsection{Pictorial representations of arc families}
\begin{figure}[ht!]
\epsfxsize = 2in
\cl{\epsfbox{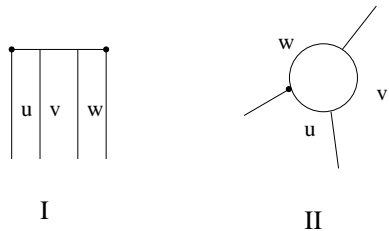}}
\caption{\label{interval}
I, bands ending on an interval,\qua II, bands ending on a circle}
\end{figure}

As explained before, there are several ways in which to imagine
weighted arc families near the boundary. 
They are illustrated in figure \ref{boundary}.
It is also convenient, to view arcs near a boundary component as coalesced into a single wide band by collapsing to a point each
interval in the window complementary to the bands; this {\it interval model} is illustrated in figure
\ref{interval}, part I.  It is also sometimes convenient to further take the image under the maps $
\ref{circleboundary}$ to produce the {\it circle model} as is depicted in figure \ref{interval}, part II.

\subsection{Gluing arc families}

\subsubsection{Gluing weighted arc families}
\label{glue}
Given two
weighted 
arc fa\-mi\-lies $(\alpha ')$ in $F_{g,m+1}^s$ and $(\beta ')$ in $F_{h,n+1}^t$
so that $\mu _i (\partial _i(\alpha '))=\mu _0(\partial _0 (\beta '))$, for some $1\leq i\leq
m$, we shall next make choices to define a weighted arc family in $F_{g+h,
m+n}^{s+t}$ as follows.  

First of all, let $\partial _0 , \partial _1, \dots , \partial _m$ denote the
boundary components of $F_{g,m+1}^s$, let  $\partial _0' , \partial _1', \dots , \partial
_n'$ denote the boundary components of $F_{h,n+1}^t$, and fix some index $1\leq i\leq m$.
Each boundary component inherits
an orientation in the standard manner from the orientations of the surfaces, and we may
choose any orientation-preserving homeomorphisms $\xi :\partial _i\to S^1$ and
$\eta :\partial _0'\to S^1$ each of which maps the initial point of the respective window to
the base-point
$0\in S^1$.  Gluing together
$\partial _i$ and $\partial _0'$ by identifying
$x\in S^1$ with $y\in S^1$ if $\xi (x)=\eta (y)$ produces a space $X$ homeomorphic to
$F_{g+h,m+n}^{s+t}$, where the two curves $\partial _i$ and $\partial _0'$ are thus
identified to a single separating curve in $X$.  There is no natural choice of
homeomorphism of $X$ with $F_{g+h,m+n}^{s+t}$, but there are canonical inclusions
$j:F_{g,m+1}^s\to X$ and $k:F_{h,n+1}^t\to X$.

We enumerate the
boundary components of $X$ in the order
$$\partial _0,\partial _1,\dots ,\partial _{i-1},\partial _1',\partial _2',\dots
\partial _{n}',\partial _{i+1},\partial _{i+2},\dots \partial _{m}$$
and re-index letting $\partial _j$, for
$j=0,1,\dots ,m+n-1$, denote the boundary components of $X$ in this order.  
Likewise, first enumerate the punctures of $F_{g,m+1}^s$ in order and then those
of $F_{h,n+1}^t$ to determine an enumeration of those of $X$, if any.  Let us choose
an orientation-preserving homeomorphism $H:X\to F_{g+h,m+n}^{s+t}$ which preserves the
labeling of the boundary components as well as those of the punctures, if any. 

In order to define the required weighted arc family, consider the partial
measured foliations ${\mathcal G}$ in $F_{g,m+1}^s$ and ${\mathcal H}$ in $F_{h,n+1}^t$
corresponding respectively to $(\alpha ')$ and $(\beta ')$.  By our assumption that $\mu _i
(\partial _i(\alpha '))=\mu _0(\partial _0 (\beta '))$, we may produce a corresponding
partial measured foliation ${\mathcal F}$ in
$X$ by identifying the points $x\in\partial _i(\alpha ')$ and
$y\in\partial _0(\beta ')$ if $c_i^{(\alpha )}(x)=c_0^{(\beta )}(y)$. 
The resulting partial measured foliation ${\mathcal F}$ may have simple closed curve leaves
which we must simply discard to produce yet another partial measured foliation
${\mathcal F}'$ in
$X$.  The leaves of ${\mathcal F}'$ thus run between boundary
components of $X$ and therefore, as in the previous section, decompose into a
collection of bands
$B_i$ of some widths $w_i$, for $i=1,2,\dots I$.  Choose a leaf of
${\mathcal F}'$ in each such band $B_i$ and associate to it the weight $w_i$ given by the width of
$B_i$ to determine a weighted arc family $(\delta ')$ in $X$ which is evidently exhaustive.  Let
$(\gamma ')=H(\delta ')$ denote the image in
$F_{g+h,m+n}^{s+t}$ under
$H$ of this weighted arc family. 

\subsubth{Lemma}\sl The $PMC(F_{g+h,m+n}^{s+t})$-orbit of
$(\gamma ')$ is well-defined as $(\alpha ')$ varies over a $PMC(F_{g,m+1}^s)$-orbit of
weighted arc families in $F_{g,m+1}^s$ and $(\beta ')$ varies over a $PMC(F_{h,n+1}^t)$-orbit
of weighted arc families in $F_{h,n+1}^t$.\rm

\proof
Suppose we are given weighted arc families $(\alpha '_2)=\phi (\alpha '_1)$,
for $\phi\in PMC(F_{g,m+1}^s)$, and
$(\beta '_2)=\psi (\beta '_1)$, for $\psi\in PMC(F_{h,n+1}^t)$, as well as 
a pair 
$H_\ell :X_\ell \to F_{g+h,m+n}^{s+t}$ of homeomorphisms as above
together with the pairs
$j_1,j_2:F_{g,m+1}^s\to X_\ell$ and
$k_1,k_2:F_{h,n+1}^t\to X_\ell$ of induced inclusions, for $\ell=1,2$.  Let
${\cal F}_\ell , {\cal F}'_\ell$ denote the partial measured foliations
and let $(\delta '_\ell)$ and $(\gamma '_\ell )$ denote
the corresponding weighted arc families in $X_\ell$ and $F_{g+h,m+n}^{s+t}$, respectively, constructed
as above from
$(\alpha '_\ell )$ and $(\beta '_\ell )$, for
$\ell =1,2$.

Let $c_\ell =j_\ell (\partial _0)=k_\ell (\partial '_i)\subseteq X_\ell$,
and remove a tubular neighborhood $U_\ell$ of $c_\ell$ in $X_\ell$ to
obtain the subsurface $X'_\ell =X_\ell -U_\ell$, for $\ell =1,2$.  Isotope 
$j_\ell ,k_\ell $ off of $U_\ell$ in the natural way to produce inclusions $j_\ell
':F_{g,m+1}^s\to X_\ell '$ and
$k_\ell ':F_{h,n+1}^t\to X'_\ell$ with disjoint images, for $\ell
=1,2$.  

$\phi$ induces a homeomorphism $\Phi :X_1'\to X_2'$ supported on
$j_1'(F_{g,m+1}^s)$ so that $j_2'\circ\phi =\Phi\circ j_1'$, and $\psi$ induces a
homeomorphism 
$\Psi :X_1'\to X_2'$ supported on $k_1'(F_{h,n+1}^t)$ so that $k_2'\circ\psi =\Psi\circ
k'_1$.  Because of their disjoint supports, $\Phi$ and $\Psi$ combine to give a homeomorphism
$G':X_1'\to X_2'$ so that $j_2'\circ \phi =G'\circ j'_1$ and $k_2'\circ\psi =G'\circ k'_1$.  We may
extend
$G'$ by any suitable homeomorphism $U_1\to U_2$ to produce a homeomorphism
$G:X_1\to X_2$.

By construction and after a suitable isotopy, $G$ maps ${\cal F} _1\cap X_1'$ to ${\cal F}_2\cap
X_2'$, and there is a power $\tau$ of a Dehn twist along $c_2$ supported on the interior of $U_2$ so
that
$K=\tau\circ G$ also maps ${\cal F}_1\cap U_1$ to ${\cal F}_2\cap
U_2$.  $K$ thus maps ${\cal F}_1'$ to ${\cal F}_2'$ and hence $(\delta _1')$ to $(\delta
_2')$.  It follows that the homeomorphism 
$$H_2\circ K\circ H_1^{-1}:F_{g+h,m+n}^{s+t}\to F_{g+h,m+n}^{s+t}$$ maps $(\gamma _1'
)$ to $(\gamma _2')$, so  $(\gamma _1')$ and $(\gamma _2')$ are indeed
in the same \hfill\break $PMC(F_{g+h,m+n}^{s+t})$-orbit.\endproof

\subsubth{Remark}
\label{rem1}
It is worth emphasizing again
that, owing to the dependence upon weights, the arcs in $\gamma '$ are {\sl not}
simply determined just from the arcs in $\alpha '$ and $\beta '$; the arcs in $\gamma '$ depend upon the weights. 
It is also worth pointing out that the composition just described is {\sl not} 
well-defined on
projectively weighted arc families but only on pure mapping class orbits of such.  In fact, by making
choices of standard models for surfaces as well as standard inclusions of these standard models,
one {\sl can} lift the composition to the level of projectively weighted arc families.  

\subsubth{Remark}
\label{rem2}
We simply discard simple closed curve components which may arise in 
our construction, and J\,L Loday
has proposed including them here in analogy to \cite{J}.  
In fact, they naturally give rise to the
conjugacy class of an element of the real group ring
 of the fundamental group of $F$.

\subsubth{Definition} \label{def1}Given $[\alpha ]\in Arc_g^s(m)$ and 
$[\beta ]\in Arc_h^t(n)$ and an index
$1\leq i\leq m$, let us choose respective deprojectivizations 
$(\alpha ')$ and $(\beta ')$ and write the weights
$$\aligned 
w(\alpha ')&=(u_0,u_1,\dots ,u_m),\\
w(\beta ')&=(v_0,v_1,\dots ,v_n).
\endaligned$$
Define
$$\aligned 
\rho _0&=\sum_{\{ b\in\beta: \partial b\cap\partial _0\neq \emptyset\}} v_i,\\
\rho _i&=\sum_{\{ a\in\alpha: \partial a\cap\partial _i\neq \emptyset\}} u_i,
\endaligned$$
where in each sum the weights are taken with multiplicity, e.g., if $a$ has both endpoints at $\partial _0$, then there are two
corresponding terms in $\rho _0$.

Since both arc families are exhaustive, $\rho _i\neq 0\neq \rho _0$, and we may re-scale 
$$\aligned
\rho _0 w(\alpha ')&=(\rho _0u_0,\rho _0u_1,\dots ,\rho _0u_m),\\
\rho _i w(\beta ')&=(\rho _iv_0,\rho _iv_1,\dots
,\rho _iv_n),\endaligned$$
so that the $0^{\rm th}$ entry of $\rho _iw(\beta ')$ agrees with the $i^{\rm th}$ entry of
$\rho _0w(\alpha ')$.  

Thus, we may
apply the composition of \ref{glue} to the re-scaled arc families to produce a
corresponding weighted arc family $(\gamma ')$ in $F_{g+h,m+n}^{s+t}$, whose projective class is
denoted
$[\gamma ]\in Arc_{g+h}^{s+t}(m+n-1)$.  We let
$$[\alpha ]\circ _i [\beta ]=[\gamma ],$$
in order to define the composition
$$\circ _i:Arc_g^s(m)\times Arc_h^t(n)\to Arc_{g+h}^{s+t}(m+n-1),
~{\rm for~any}~i=1,2,\dots ,m.$$ As in Remark \ref{rem1}, this
composition is {\sl not} a simplicial map but just a topological one.

\subsubsection{A pictorial representation of the gluing}
\begin{figure}[ht!]
\epsfxsize = \textwidth
\epsfbox{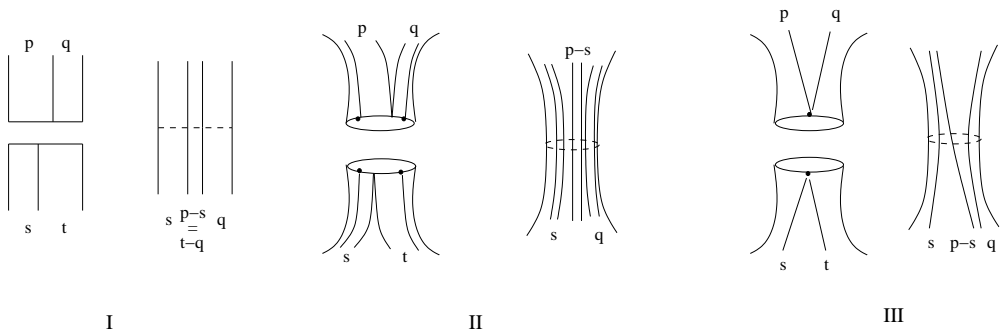}

\caption{\label{gluefig}
The gluing: I, in the interval picture,  II, in the windows with bands
picture and III, in the arcs running to a marked point version.}
\end{figure}

A graphical representation of the gluing can be found in 
figure \ref{gluefig},  where we present the gluing in three of the different 
models. More examples of gluing in the interval and circle
models  can be found in figure \ref{bands} and throughout the later
sections.

\subsection{The basic topological operads}

\subsubth{Definition} For each $n\geq 0$, let $\Arc_{cp}(n)=Arc_0^0(n)$ (where the ``$cp$'' stands for
compact planar), and furthermore, define the direct limit
$\Arc(n)$ of $Arc_g^s(n)$ as $g,s\rightarrow \infty$ under the natural inclusions.

\subsubth{Theorem} \sl The compositions $\circ_{i}$ of Definition \ref{def1} imbue
the collection of spaces $\Arc_{cp}(n)$ with the structure of a
topological operad under the natural $\mathbb{S}_{n}$--action on labels on the boundary components.
The operad has a unit $1\in \Arc_{cp}(1)$ given by the class of an arc in the cylinder meeting both
boundary components, and the operad is cyclic for the natural
$\mathbb{S}_{n+1}$--action. \rm

\proof The first statement follows from the standard 
operadic manipulations. The second statement is immediate from the definition of composition.
For the third statement, notice that the composition treats the two surfaces
symmetrically, and the axiom for cyclicity again follows from the standard
operadic manipulations.\endproof

In precisely the same way, we have the following theorem.

\subsubth{Theorem}\sl \label{opinfty}The composition $\circ _i$ of Definition \ref{def1}
induces a composition
$\circ _i$ on $\Arc(n)$ which imbues this collection of spaces with the natural
structure of a cyclic topological operad with unit.\rm

\subsubsection{The deprojectivized spaces $\darc$}
For the following it is convenient to introduce deprojectivized
arc families. This amounts to adding a factor ${\mathbb R}_{> 0}$
for the overall scale.

Let $\darc(n)=\Arc(n) \times {\mathbb R}_{> 0}$ be the space
of weighted arc families;  
it is clear that $\darc(n)$ is homotopy equivalent to $\Arc(n)$.

As the  definition \ref{def1} of 
gluing was obtained by lifting to weighted arc families
and then projecting back, we can promote the compositions to the level
of the spaces $\darc(n)$.
This endows the spaces $\darc(n)$ 
with a structure of a {cyclic} operad as well. Moreover, by
construction the two operadic structures are compatible.
 This type of composition 
can be compared to the composition of loops, where such a rescaling is also
inherent. In our case, however, 
the scaling is performed on both sides which renders the operad cyclic.

In this context, the total weight at a given boundary
component given by the sum of the individual weights $w_t$ of 
incident arcs makes sense,
and thus the map \ref{circleboundary} can be naturally viewed as a
map to a circle of radius $\sum_t w_t$.

\subth{Notation} 
We denote the operad on the collection of spaces \hfill\break
$\Arc(n)$ by $\Arc$ and the
operad on the collection of spaces $\darc(n)$ by $\darc$.  By an ``Arc
algebra'', we mean an algebra over the homology operad of $\Arc$.
Likewise, $\Arc _{cp}$ and $\darc _{cp}$ are comprised of the spaces
$\Arc_{cp}(n)$ and $\darc_{cp}(n)$ respectively, and an ``Arc$_{cp}$ algebra'' is an algebra
over the homology of $\Arc _{cp}$

\subsubth{Remark}
\label{rem3}
In fact, the restriction that the arc families under consideration must be
exhaustive can be relaxed in several ways with the identical definition 
retained for the composition.
For instance, if the projectively weighted arc family 
$[\alpha ]\in Arc _g^s(m)$ fails to meet
the $i^{\rm th}$ boundary component and $[\beta ]\in Arc_h^t(n)$, then 
we can set $[\alpha ]\circ _i
[\beta ]$ to be $[\alpha  ]$ regarded as an arc family in
$F_{g,m+1}^s\subseteq F_{g+h,m+n}^{s+t}$, where the boundary 
components have been
re-labeled. 

In both formulations, 
we must require that $\partial _0 [\beta ]\neq\emptyset$ in
order that $[\alpha ]\circ _i
[\beta ]$ is non-empty, but this 
asymmetric treatment destroys cyclicity.  

Another possibility is to
include the empty arc family in the operad as any arc family with all
weights zero to preserve cyclicity.   The composition then imbues the
deprojectivized $Arc(F_{g,n+1}^s)$ themselves 
with the
structure of a topological operad, 
but the corresponding homology operad, in light of the Sphericity
Conjecture \ref{sphere} would 
be trivial i.e.,\ the trivial one-dimensional $\Sn$--module
for each $n$.

\subsubsection{Turning on punctures and genus}

\begin{figure}[ht!]
\epsfxsize = 4in
\cl{\epsfbox{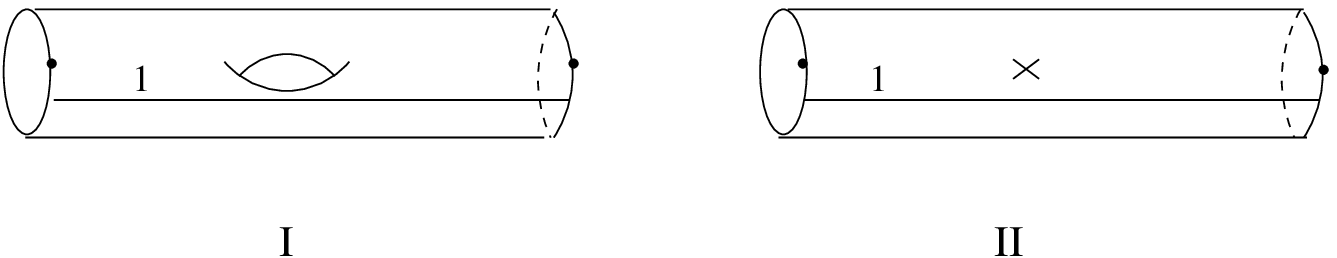}}
\caption{\label{topops} I, the genus generator,\qua II, the puncture operator}
\end{figure}

When allowing genus or the number of punctures to be different
from zero, there are two operators for the topological
operad that generate arc families of all genera respectively with any 
number of punctures from arc families of genus zero with no punctures.
These operators are depicted in figure \ref{topops}.

For the corresponding linear homology operad, we expect similarities with
\cite{LM}, in which the punctures play the role of the second
set of points. This means that the linear operad will act on tensor
powers of two linear spaces, one for the boundary components and one for
the punctures. In this extension of the operadic framework, 
there is no gluing on the punctures and
the respective linear spaces should be regarded as parameterizing 
deformations which
perturb each operation separately.

\subsubsection{Other Props and Operads}
\label{rem4}
There are related sub-operads and sub-props of $Arc(F_{g,n+1}^s)$ other 
than $\Arc_{cp}$
and $\Arc$
which are of interest.  In the general case, one may specify a symmetric 
$(n+1)$-by-$(n+1)$ matrix
$A^{(n)}$  as well as an $(n+1)$-vector $R^{(n)}$ of zeroes and ones 
over $\ZZ$ and consider the subspace of
$Arc(F_{g,n+1}^s)$ where arcs are allowed to run between boundary 
components $i$ and $j$
if and only if
$A_{ij}^{(n)}\neq 0$ and are required to meet boundary component $k$ if 
and only if $R_k^{(n)}\neq
0$.  For instance, the case of interest in this paper corresponds to
$A^{(n)}$ the matrix and $R^{(n)}$ the vector 
whose entries are all one.  In Remark \ref{rem3}, we also mentioned the
example with $A^{(n)}$ the matrix consisting entirely of entries one and
$R^{(n)}$ the standard first unit basis vector.

For a class of examples, consider a partition of
$\{ 0,1,\dots n\}  =I^{(n)}\sqcup O^{(n)}$, into 
``inputs'' and ``outputs'', where $A_{ij}^{(n)}=1$
if and only if
$\{ i,j\}\cap I^{(n)}$ and
$\{ i,j\}\cap O^{(n)}$ are each singletons, 
and $R^{(n)}$ is the vector whose entries are all one.
The corresponding prop is presumably related to the string prop of \cite{CS}.

There
are many other interesting possibilities.  For instance,
the cacti operad of \cite{V} and the spineless cacti of \cite{K}
will be discussed and studied in sections 2 and 3, and other related
cactus-like examples are studied in \cite{K}.
A further variation is to stratify the space by
incidence matrices with non--negative integer entries corresponding
to the number of arcs at each boundary component.

\section{The Gerstenhaber and BV structure of $\Arc$}
\label{BV}

In this section we will show that there is a structure of Batalin-Vilkovisky algebra on
Arc algebras and Arc$_{cp}$ algebras. More precisely, any algebra over
the singular chain complex operad $C_*(\Arc _{cp})$ or 
$C_*(\Arc )$ is a  BV algebra up to certain chain homotopies. For convenience, $C_*(n)$ will denote
$C_*(\Arc _{cp}(n))$ or $C_*(\Arc(n))$, since  for our calculations, we only 
require $\Arc _{cp}$. An element in $C_*$ will be also called an arc family.
We shall realize the underlying surfaces $F_{0,n+1}^0$ in the plane with the
boundary $0$ being the outside circle.
The conventions we use for the drawings is that all circles are
oriented counterclockwise. Furthermore if there is only one
marked point on the circle, it is the beginning of the window. If there are two 
marked points on the boundary, the window is the smaller arc between the two.

Via gluing, any arc family in $C_*(n)$ gives rise to an n--ary
operation on arc families in $C_*$. Here one has to be careful with the
parameterizations.  To be completely explicit we will always include
them if we use a particular family as an operation.

As mentioned in the Introduction, many of the calculations of this section are inspired by \cite{CS}.

\subsection{A reminder on some algebraic structures}

In this section, we would like to recall some basic definitions 
of algebras and their relations which we will employ in the following.
The proofs in this subsection 
are omitted, since they are straightforward computations.
They can be found in \cite{G,Ge}.

\subsubth{Definition}(Gerstenhaber)\qua
{A pre--Lie algebra is a $\Zz$ graded vector space $V$ together with a bilinear
operation $*$ that satisfies
\begin{equation*}
(x*y)*z-x*(y*z) = (-1)^{|y||z|} ((x*z)*y -x*(z*y))
\end{equation*}}
Here $|x|$ denotes the $\Zz$ degree of $x$.

\subsubth{Definition}
{An odd Lie algebra is a  $\Zz$ graded vector space $V$ together with
a bilinear operation $\{\, ,\, \}$ which satisfies
$$\{a \, ,\,  b\} = (-1)^{(|a|+1)(|b|+1)}\{b \, ,\,  a\}\leqno{(1)}$$
$$\{a \, ,\,  \{b \, ,\,  c\} \} = 
\{\{a \, ,\,  b\} \, ,\,  c \}
+(-1)^{(|a|+1)(|b|+1)} \{b \, ,\,  \{a \, ,\,  c\} \}\leqno{(2)}$$
}

\subsubth{Remark} Pre--Lie algebras of the above type are sometimes
also called right symmetric algebras. There is also the notion 
of a left symmetric algebra, which satisfies
\begin{equation*}
(x*y)*z-x*(y*z) = \pm (y*x)*z- y*(x*z)
\end{equation*}
Given a left symmetric algebra its opposite algebra is right symmetric and
vice--versa. Here the multiplication for the opposite algebra
$A^{opp}$ is $a*^{opp}b= b*a$. Hence it is a matter of taste, which
algebra type one regards. To match with string topology, one has to use
left symmetric algebras, while to match with the Hochschild cochains,
one will have to  use right symmetric algebras.

\subsubth{Definition} {A Gerstenhaber algebra or an odd Poisson
algebra is a $\Zz$ graded, graded--commutative associative
algebra  $(A,\cdot)$ endowed with an odd Lie algebra structure 
$\{ \, ,\, \}$ which satisfies the compatibility equation
$$
\{a \, ,\,  b \cdot c \}= \{a \, ,\,  b\} \cdot c+
(-1)^{|b|(|a|+1)} b\cdot \{a \, ,\,  c\} 
$$
}

\subsubth{Proposition}
(Gerstenhaber)\qua
{\sl For any pre-Lie algebra $V$, the 
bracket $\{\, \, ,\,  \,\}$ defined by 
\begin{equation}
\label{bracketdef}
\{a \, ,\,  b\}:= a*b - (-1)^{(|a|+1)(|b|+1)}b*a
\end{equation}
endows $V$ with a structure 
of odd Lie algebra.}

\subsubth{Remark} The same holds true for a right symmetric algebra.

\subsubth{Definition}
{A Batalin--Vilkovisky (BV) algebra is an associative super--commutative
algebra $A$ together with an operator 
$\Delta$ of degree $1$ that satisfies
\begin{eqnarray*}
\Delta^2&=&0\\
\Delta(abc) &=&\Delta(ab)c+(-1)^{|a|}  a\Delta(bc) + (-1)^{|sa||b|} 
b\Delta (ac)\\
& -&\Delta(a)bc
-(-1)^{|a|} a\Delta(b)c-(-1)^{|a|+|b|}
ab\Delta(c)
\end{eqnarray*}
}	

Here super--commutative means as usual $\Zz$ graded commutative,
i.e.\ $ab= (-1)^{|a||b|}ba$.

\subsubth{Proposition}(Getzler)\qua
{\sl For any BV--algebra $(A,\Delta)$ define
\begin{equation}
\label{gbv}
\{a \, ,\,  b\}:= (-1)^{|a|} \Delta(ab)-(-1)^{|a|} \Delta(a)b-a \Delta(b)
\end{equation}
Then $(A, \{\, \, ,\,  \,\})$ is a Gerstenhaber algebra.}

\subsubth{Definition}
We call a triple $(A, \{\, \, ,\,  \,\}, \Delta)$ a GBV--algebra
if $(A, \Delta)$ is a BV algebra and 
$ \{\, \, ,\,  \,\}:A\otimes A \rightarrow A$ 
 satisfies
the equation (\ref{gbv}). By the Proposition above $(A,\{\, \, ,\,  \,\})$
is a Gerstenhaber algebra. \medskip

The purpose of this definition is that in some cases as in the case of
the Arc operad it happens that a bracket as well as the BV operator
appear naturally. In our case the bracket comes naturally from a
pre--Lie structure which one can also view as a $\cup_1$ product,
while the BV operator appears from a cycle naturally parameterized by
$S^1=Arc(F_{0,2}^0=\Arc_{cp}(1)$.  We use the name GBV algebra to
indicate that these structures although having independent origin are
indeed compatible. As we will show below, the independent origin can
be interpreted as saying that the Gerstenhaber structure is governed
by one suboperad and the BV--structure by another bigger suboperad
which contains the previous one. Moreover the bracket of the BV
structure coincides with the bracket of the Gerstenhaber structure
already present in the smaller suboperad.

\subsection{Arc families and their induced operations}

\begin{figure}[ht!]
\epsfxsize =\textwidth
\epsfbox{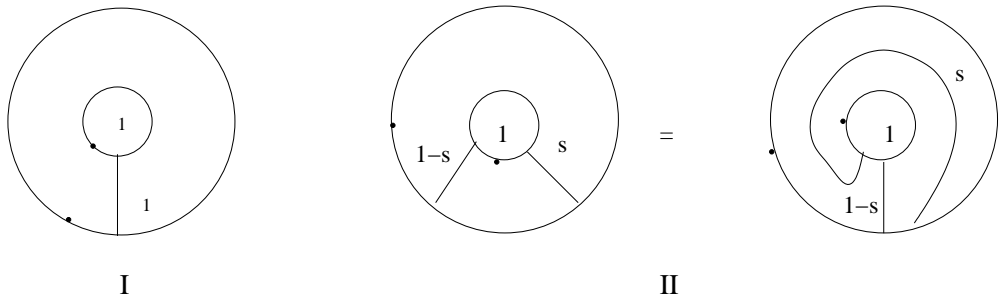}
\caption{\label{1ops}
I, the identity and II, the arc family $\delta$ yielding the BV operator}
\end{figure}

The points in $\Arc _{cp}(1)$ are parameterized by the circle, which is identified with $[0,1]$, where $0$ is
identified to $1$. To describe a parameterized family of weighted arcs, we shall specify weights that depend
upon the parameter $s\in [0,1]$. Thus, by taking $s\in [0,1]$ figure \ref{1ops} describes a cycle 
$\delta\in C_1(1)$
that spans $H_1(\Arc _{cp}(1))$.

\begin{figure}[ht!]
\epsfxsize = \textwidth
\epsfbox{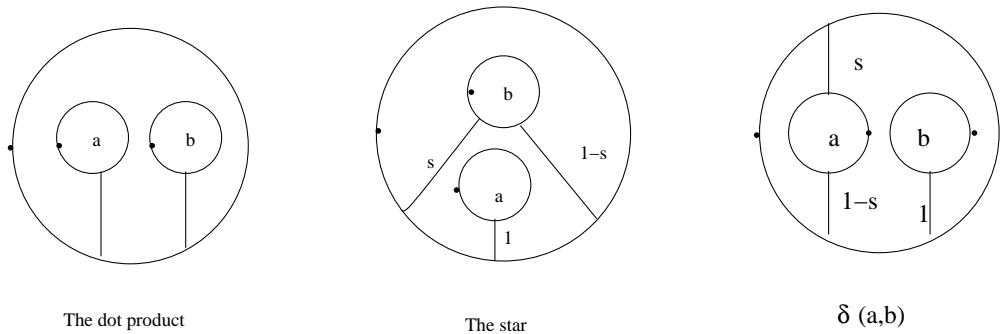}
\caption{\label{2ops}The binary operations}
\end{figure}

As stated above, there is an operation associated to the family $\delta$.
For instance, if $F_1$ is any arc family $F_1:k_1 \rightarrow \Arc _{cp}$,
$\delta F_1$ is the family parameterized by 
$I \times k_1 \rightarrow \Arc _{cp}$ with the map given by the picture by inserting
$F_1$ into the position 1. By definition, $$\Delta=-\delta \in C_1(1).$$
In $C_*(2)$ we have the basic families depicted in
figure \ref{2ops} which in turn yield operations
on $C_*$.

To fix the signs, we fix the parameterizations we 
will use for the glued families as follows:
say the families $F_1,F_2$ are parameterized by $F_1:k_1 \rightarrow \Arc _{cp}$
and $F_2:k_2 \rightarrow\Arc _{cp}$ and $I=[0,1]$.
Then $F_1\cdot F_2$ is the family parameterized by 
$k_1 \times k_2 \rightarrow \Arc _{cp}$ as defined by figure \ref{2ops} (i.e.,
the arc family $F_1$ inserted in boundary $a$ and the arc family
$F_2$ inserted in boundary $b$).

\begin{figure}[ht!]
\epsfxsize = 4in
\cl{\epsfbox{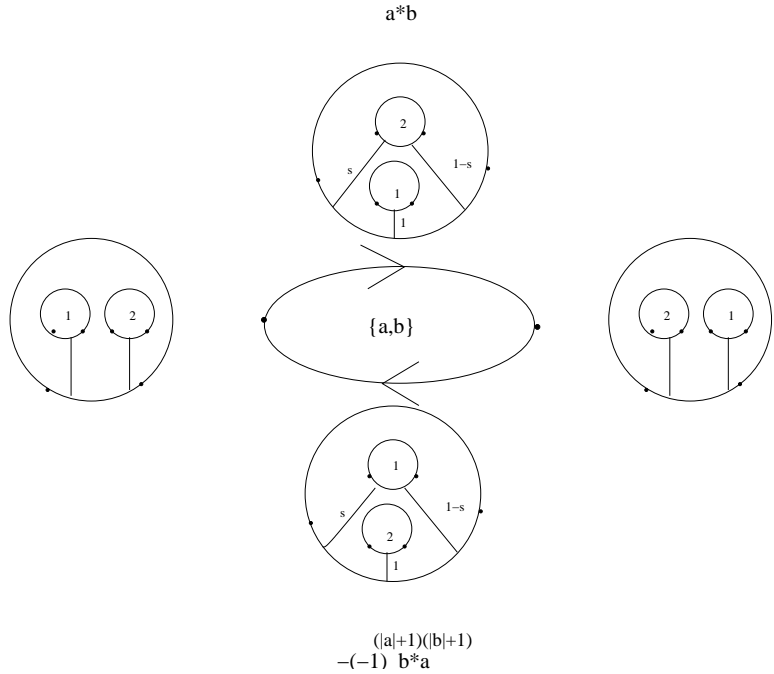}}
 \caption{\label{bracket}The definition of the Gerstenhaber bracket}
\end{figure}

Interchanging labels 1 and 2 and using $*$ as a chain homotopy as in figure \ref{bracket} 
yields the commutativity of $\cdot$ up to chain homotopy
\begin{equation}
d (F_1*F_2)=(-1)^{|F_1||F_2|}F_2\cdot F_1-F_1\cdot F_2
\end{equation}
Notice that the product $\cdot$  is associative up to chain homotopy.

Likewise $F_1 * F_2$ is defined to be the operation given by the second
family of figure \ref{2ops} with $s \in I = [0,1]$ parameterized over
$k_1 \times I \times k_2 \rightarrow \Arc _{cp}$.

By interchanging the labels, we can produce a cycle $\{F_1, F_2\}$ 
 as shown in figure \ref{bracket}
where now the whole family is parameterized by 
$k_1 \times I \times k_2 \rightarrow \Arc _{cp}$.
$$\{F_1,F_2\} := F_1*F_2 - (-1)^{(|F_1|+1)(|F_2|+1)}F_2*F_1.$$

\subsubth{Definition} We have
defined the following elements in $C_*$: 

$\delta$ and $\Delta=-\delta$ in  $C_1(1)$;

$\cdot$ in $C_0(2)$, which is commutative and associative up to a boundary. 

$*$ and $\{-,-\}$ in $C_1(2)$ with $d(*)=\tau \cdot - \cdot$ and
$\{-,-\}=*-\tau *$.

Note that $\delta, \cdot$ and $\{-,-\}$ are cycles, whereas $*$ is not.

\subsubth{Remark} We would like to point out that the symbol 
$\bullet$ in the
standard super notation of odd Lie brackets $\{ a \bullet b \}$, which
is assigned to have 
an intrinsic degree of 1, corresponds geometrically in our situation
to the one--dimensional interval $I$.

\subsection{The BV operator}

The operation corresponding to the arc family $\delta$ is easily seen
to square to zero in homology. It is therefore a differential
and a natural candidate for a derivation or a 
higher order differential operator.
It is easily checked that it is {\em not a derivation}, but it {\em is a
BV operator}, as we shall demonstrate. 

\begin{figure}[ht!]
\epsfxsize =4.5in
\cl{\epsfbox{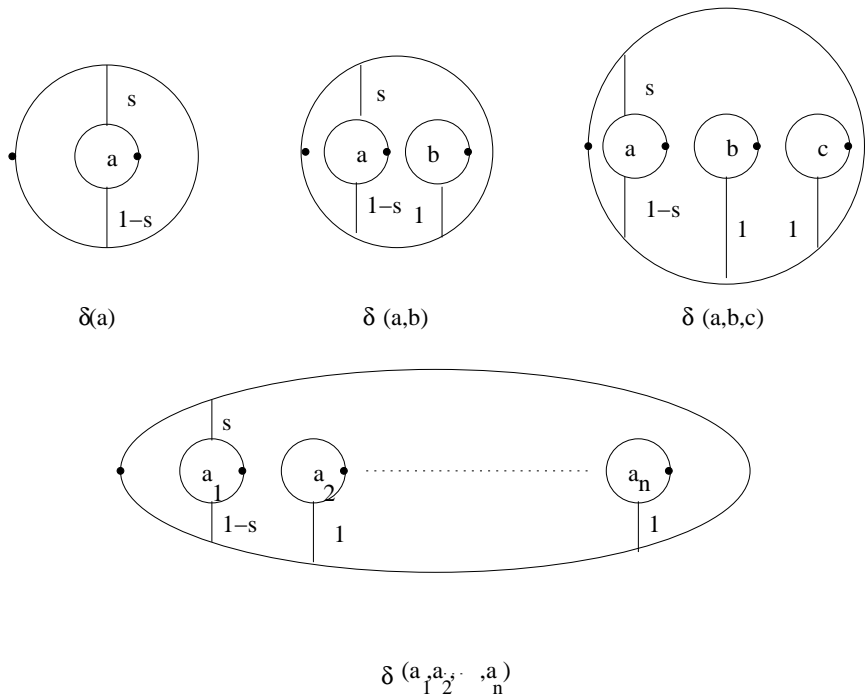}}
\caption{\label{delfig}The definition of the $n$--ary operations $\delta$}
\end{figure}

It is convenient to introduce the family of operations on arc families
$\vardel$ which are defined by figure \ref{delfig},
where the families are parameterized over $I \times k_{a_1} \times \dots
\times k_{a_n}$.

We notice the following relations which are the raison d'\^etre for
this definition:
\begin{eqnarray}
\label{threedel}
%\delta(a) &\sim& \vardel(a)\nn\\
\delta(ab) &\sim& \vardel(a,b) + (-1)^{|a||b|}\vardel(b,a)\nn\\
\delta(abc) &\sim& \vardel(a,b,c)+(-1)^{|a|(|b|+|c|)}\vardel(b,c,a)\\
&&+(-1)^{|c|(|a|+|b|)}\vardel(c,a,b)\nn\\
\delta(a_1 a_2\cdots a_n)&\sim& \sum_{i=0}^{n-1} (-1)^{\s(c^i,a)} 
\delta(a_{c^i(1)}, \dots, a_{c^i(n)})
\end{eqnarray}
where $c$ is the cyclic permutation $(1,\ldots ,n)$ and  
$\s(c^i,a)$ is the standard supersign of the permutation. The homotopy here
is just a reparameterization of the variable $s\in I$.

There is a further relation immediate from the definition which 
shows that the only ``new'' operation is $\vardel(a,b)$
\begin{equation}
\vardel(a_1,a_2, \dots, a_n) \sim \vardel(a_1, a_2 a_3 \cdots a_n)
\end{equation}
where we use a homotopy to scale all weights of the bands not hitting
the boundary 1 to the value 1.

\subsubth{Lemma}
\label{partbv}
\begin{equation}
\vardel(a,b,c) \sim (-1)^{(|a|+1)|b|} b \vardel(a,c)
+\vardel(a,b)c -\vardel(a)bc
\end{equation}

\begin{figure}[ht!]
\epsfxsize = 4.8in
\cl{\epsfbox{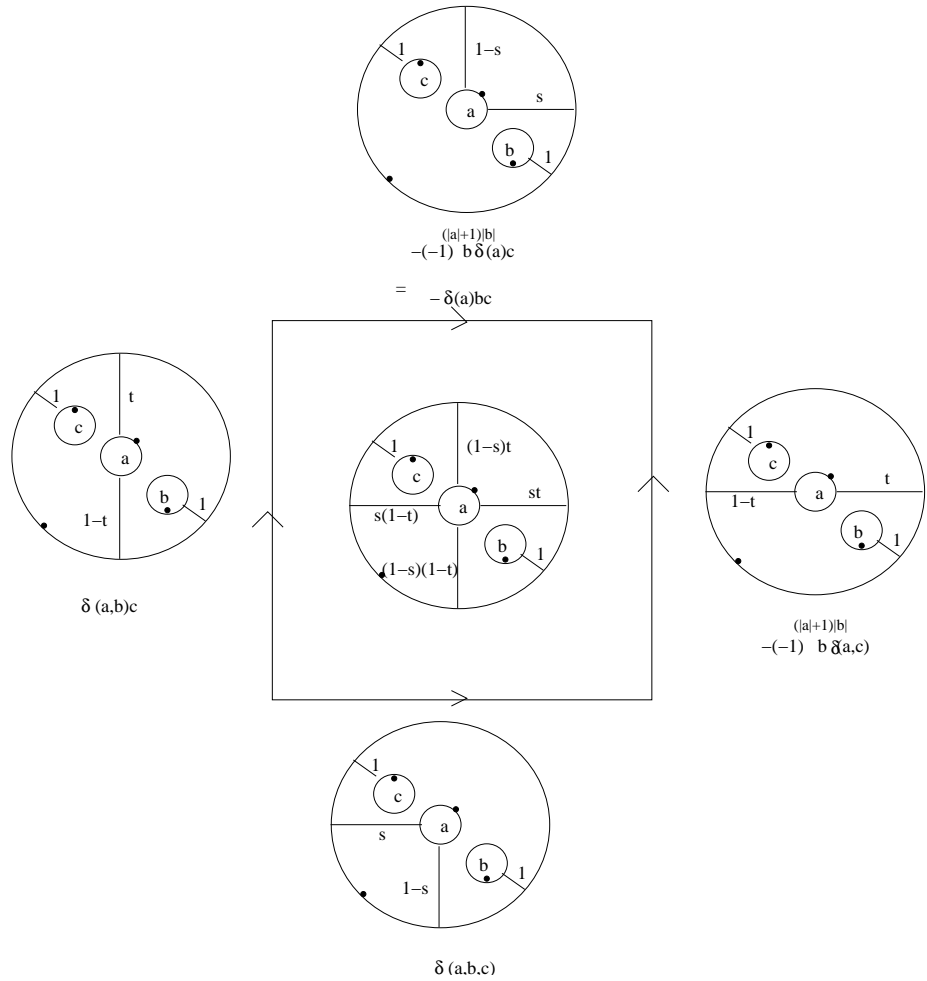}}
\caption{\label{partbvfig}The basic chain homotopy responsible for BV}
\end{figure}

\proof The proof is contained in figure \ref{partbvfig}. Let
$a:k_a \rightarrow \Arc _{cp}$,
$b:k_b \rightarrow \Arc _{cp}$ and  $c:k_c \rightarrow \Arc _{cp}$, be arc families
then the two parameter family
filling the square is parameterized over 
$I \times I \times k_a \times k_b \times k_c$. This family
gives us the desired chain homotopy.\endproof

\subsubth{Proposition}
\sl The operator $\Delta$  satisfies the relation
of a BV operator up to chain homotopy.
\begin{eqnarray}
\Delta^2&\sim&0\nn\\
\Delta(abc)&\sim&\Delta(ab)c+(-1)^{|a|}  a\Delta(bc) + (-1)^{|sa||b|} b\Delta (ac)
 -\Delta(a)bc\nn\\
&&-(-1)^{|a|} a\Delta(b)c-(-1)^{|a|+|b|}
ab\Delta(c)
\end{eqnarray}
Thus, any Arc algebra and any Arc$_{cp}$ algebra is a BV algebra.
\rm

\proof The proof follows algebraically from Lemma \ref{partbv} and 
equation (\ref{threedel}). We can also make the chain homotopy 
explicit. This has the advantage of illustrating the symmetric nature of
this relation in $C_*$ directly.

\begin{figure}[ht!]
\epsfxsize = \textwidth
\epsfbox{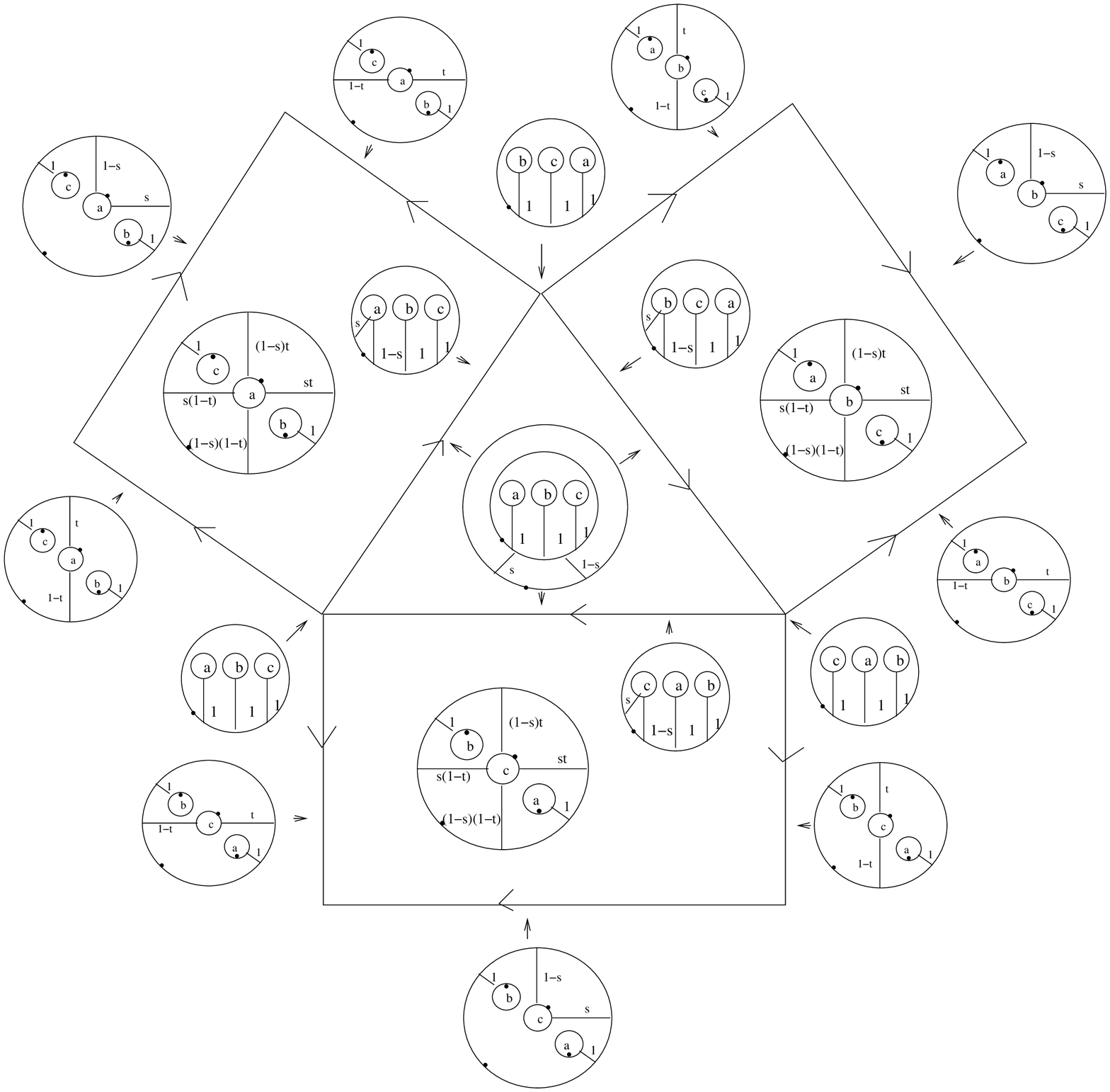}
\caption{\label{bvdirect}The homotopy BV equation}
\end{figure}

 Given arc families $a:k_a \rightarrow \Arc _{cp}$,
 $b:k_b \rightarrow \Arc _{cp}$ and  $c:k_c \rightarrow \Arc _{cp}$, 
we define the two parameter family  defined by the figure \ref{bvdirect}
where the families in the rectangles are  
 the depicted two parameter families  parameterized over 
$I \times I \times k_a \times k_b \times k_c$ and
the triangle is not filled, but rather its boundary is the operation
$\delta(abc)$.

From the diagram we get the chain homotopy consisting of three, and respectively
twelve, terms. 
\begin{eqnarray}
\delta(abc) &\sim& \vardel(a,b,c)+(-1)^{|a|(|b|+|c|)}\vardel(b,c,a)
+(-1)^{|c|(|a|+|b|)}\vardel(c,b,a)\nn\\
&\sim&(-1)^{(|a|+1)|b|} b \vardel(a,c)
+\vardel(a,b)c -\vardel(a)bc
+ (-1)^{|a|} a \vardel(b,c)\nn\\
&&+(-1)^{|a||b|}\vardel(b,a)c -(-1)^{|a|}a\vardel(b)c
+(-1)^{(|a|+|b|)|c|} a \vardel(b,c)\nn\\
&&+(-1)^{|b|(|a|+1|)+|a||c|}b\vardel(c,a)c - (-1)^{|a|+|b|} ab\vardel(c)\nn\\
&\sim&\delta(ab)c+(-1)^{|a|}  a\delta(bc) + (-1)^{|a+1||b|} b\delta (ac)
 -\delta(a)bc\nn\\
&&-(-1)^{|a|} a\delta(b)c-(-1)^{|a|+|b|}ab\vardel(c)
\end{eqnarray}

\subsection{Gerstenhaber Structure}

We have already defined the 
operation whose odd commutator is given by the BV operator.

\subsubth{Theorem}\sl
The Gerstenhaber bracket induced by $\Delta$ is gi\-ven by the operation
$$
\{a,b\} = a*b - (-1)^{(|a|+1)(|b|+1)}b*a
$$
In other words $\Arc$ is a GBV--Algebra up to homotopy.
\rm

\begin{figure}[ht!]
\epsfxsize = \textwidth
\epsfbox{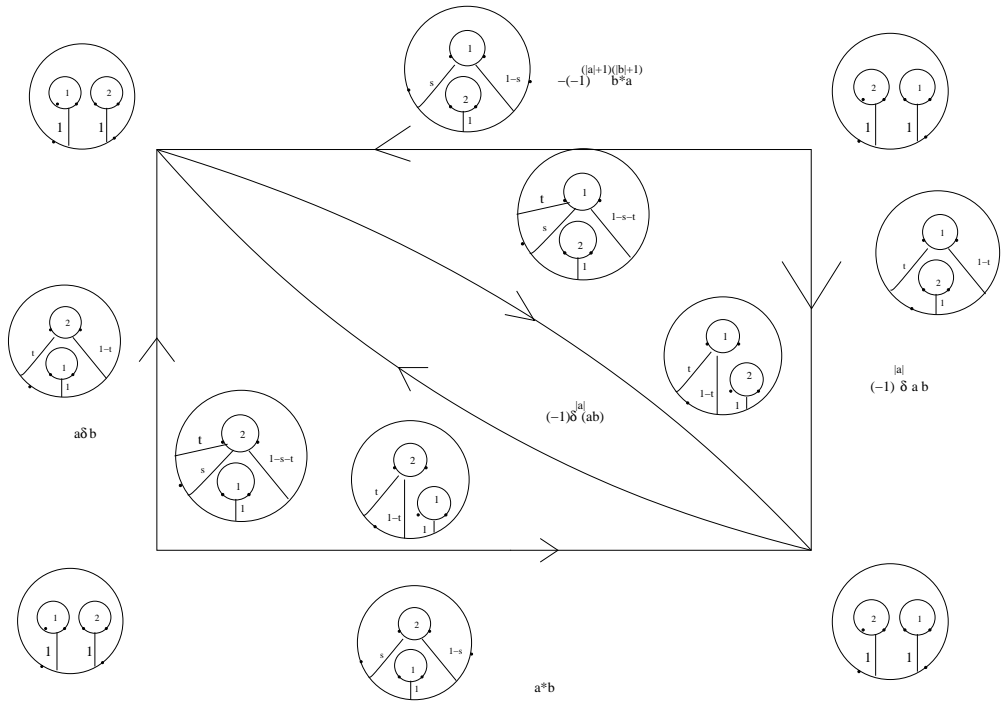}
\caption{\label{bv} The odd commutator realization of the bracket}
\end{figure}
\proof We consider the  arc family depicted  in figure \ref{bv},
where the two parameter arc families are parameterized via
$k_1 \times T \times k_2 \rightarrow \Arc _{cp}$, where 
$T$ is the triangle $T:=\{(s,t)\in [0,1]\times [0,1]:s+t\leq 1\}$ with
induced orientation,
and where the middle loop is $(-1)^{|F_1|}\delta(F_1 \cdot F_2)$ 
suitably re--parameterized. Inserting $F_1$ in the boundary 1,
and $F_2$ in boundary
2, and passing to homology, we can read off the relation:
\begin{multline*}
(-1)^{|F_1|}\delta(F_1\cdot  F_2)=\\ (-1)^{(|F_1|+1)(|F_2|+1)}F_1*F_2
+ (-1)^{|F_1|} \delta(F_1)\cdot F_2
-F_1 *F_2
+F_1\cdot (\delta F_2)
\end{multline*}
or with $\Delta=-\delta$:
$$
\{F_1,F_2\}= (-1)^{|F_1|}\Delta(F_1\cdot  F_2) 
-(-1)^{|F_1|} \Delta(F_1)\cdot F_2
-F_1\cdot (\Delta F_2)
$$

\subsubth{Remark} Algebraically, the Jacobi identity and the 
derivation property of the bracket follow from the BV relation.
For the less algebraically inclined we can again make everything
topologically explicit. This also has the virtue of showing how 
different weights contribute to topologically distinct gluings. 
This treatment also  
shows that we can restrict ourselves to the case of linear (Chinese)
trees of \S \ref{cacti} or to cacti without
spines (cf.\ \S5 and \cite{K}) and therefore have a Gerstenhaber structure on
this level.

\subsubsection{The associator}
\begin{figure}[ht!]
\epsfxsize = \textwidth
\epsfbox{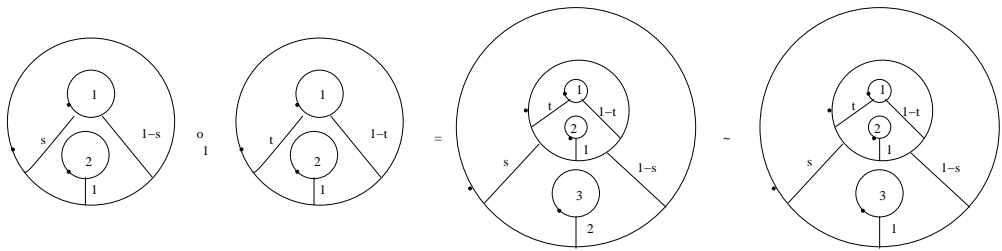}
\caption{\label{circ1} The first iterated gluing of $*$}
\end{figure}

It is instructive to do the calculation 
in the arc family picture with the operadic notation.
For the gluing $* \circ _1 *$ 
we obtain the elements in $C_2(2)$ presented in 
figure \ref{circ1} to which we apply the homotopy of changing
the weight on the boundary 3 from 2 to 1 while keeping
everything else fixed. We call this normalization.

\begin{figure}[ht!]
\epsfxsize = \textwidth
\epsfbox{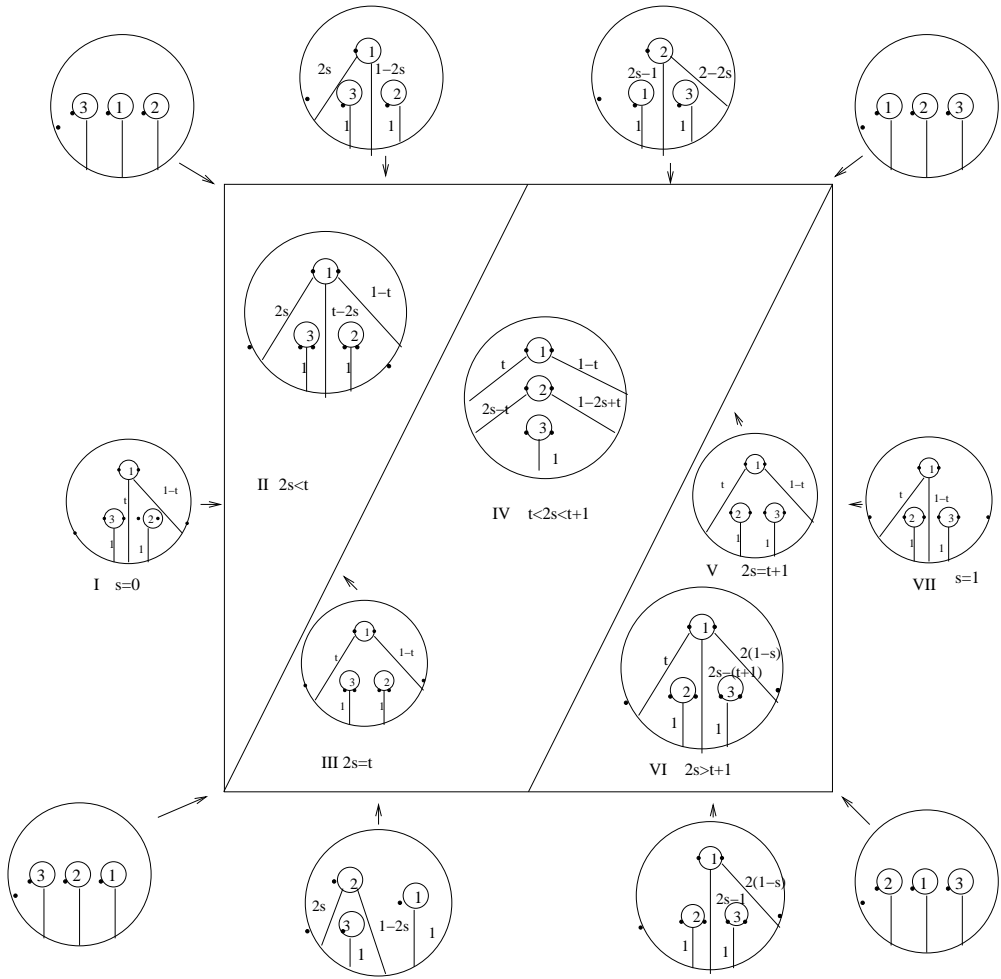}
\caption{\label{circ1b} The glued family after normalization}
\end{figure}

\begin{figure}[ht!]
\epsfxsize = \textwidth
\epsfbox{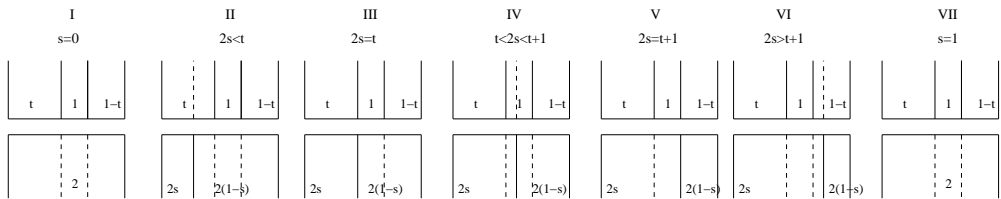}
\caption{\label{bands}The different cases of gluing the bands}
\end{figure}

Unraveling the definitions for the normalized version yields
figure \ref{circ1b},
where in the different cases the gluing of the bands is shown in
figure \ref{bands}.

\begin{figure}[htb]
\epsfxsize = \textwidth
\epsfbox{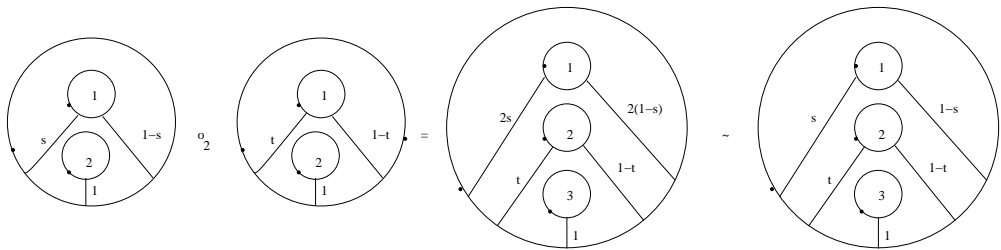}
\caption{\label{circ2}The other iteration of $*$}
\end{figure}

The gluing $* \circ_2 *$ in arc families is simpler and yields the
gluing depicted in figure \ref{circ2}  to which we apply a normalizing 
homotopy --- 
by changing the weights on the bands emanating from boundary 
1 from the pair $(2s,2(1-s))$ to $(s,(1-s))$ using pointwise the  
homotopy $(\frac{1+t}{2}2s,\frac{1+t}{2}(1-s))$ for $t\in [0,1]$:

Combining figures \ref{circ1b} and \ref{circ2} 
while keeping in mind the parameterizations we
can read off the pre--Lie relation:
\begin{multline}
F_1 * (F_2 *F_3) - (F_1*F_2)*F_3\sim \\
(-1)^{(|F_1|+1)(|F_2|+1)}(F_2*(F_1*F_3)- (F_2*F_1)*F_3)
\end{multline}
which shows that the associator is symmetric in the first two variables and
thus following Gerstenhaber \cite{G} we obtain:
\subsubth{Corollary}\sl
\label{jac} The bracket
$\{\, , \, \}$
satisfies the odd Jacobi identity.\rm

\subsubsection{The Gerstenhaber structure}

\begin{figure}[ht!]
\epsfxsize = \textwidth
\epsfbox{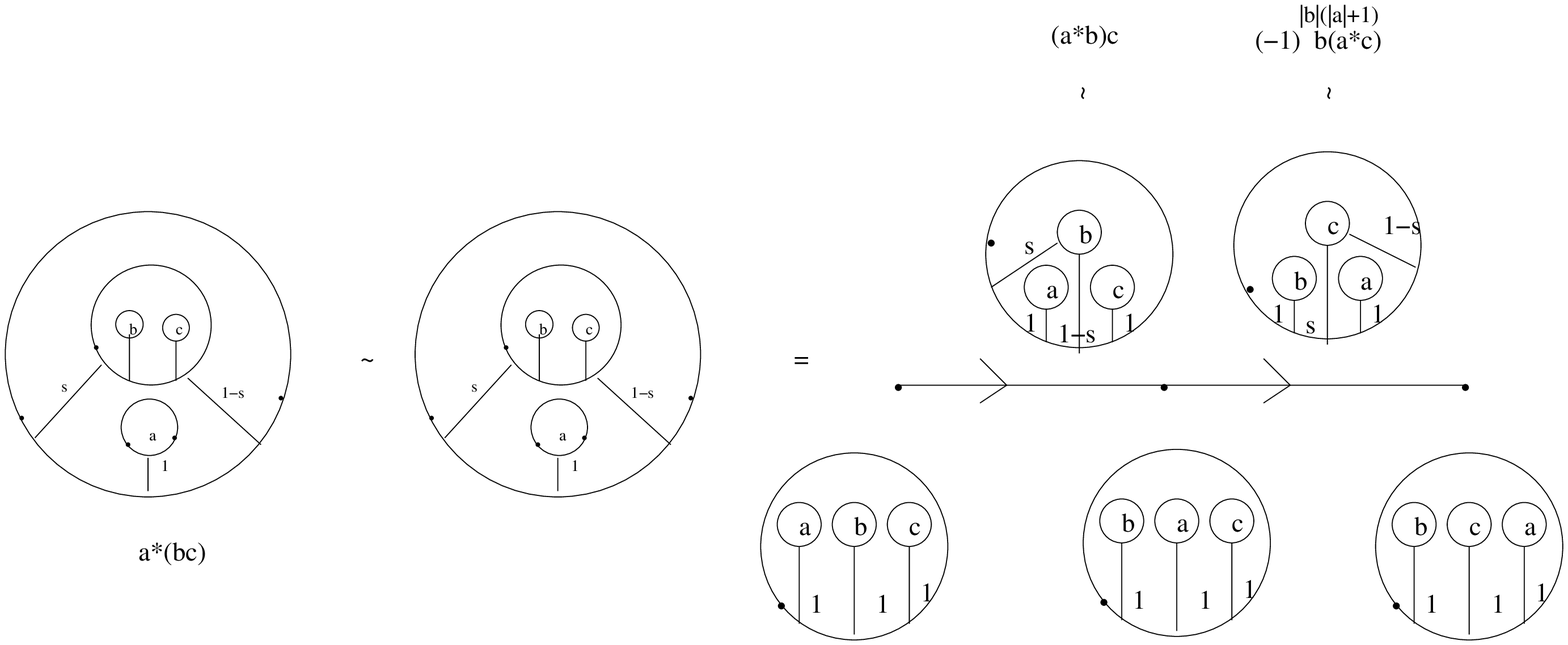}
\caption{\label{compat} The first compatibility equation}
\end{figure}

\begin{figure}[ht!]
\epsfxsize = 4.5in
\cl{\epsfbox{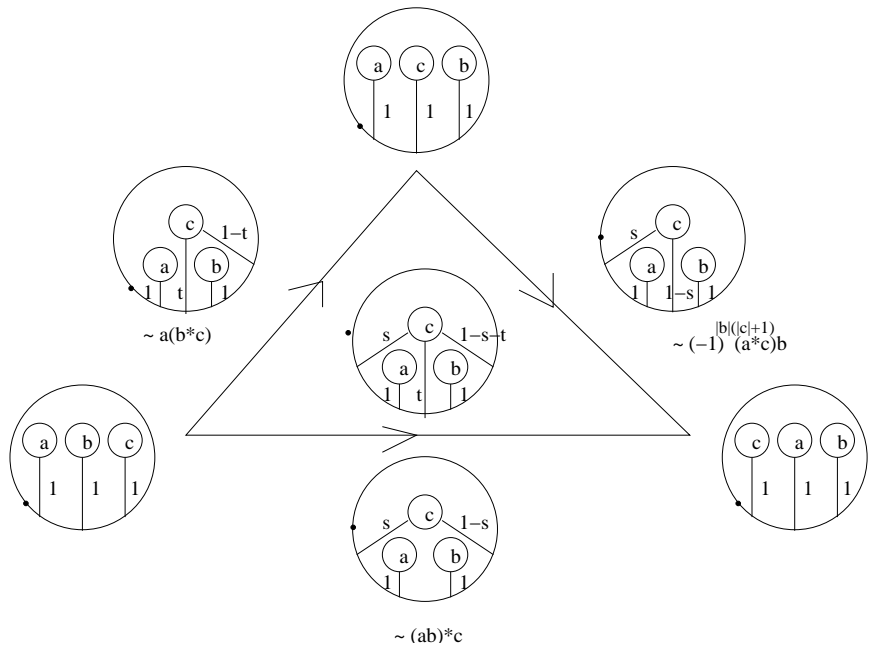}}
\caption{\label{compat2} The second compatibility equation}
\end{figure}

The derivation property of the bra\-cket follows from 
the compatibility equations which are proved by the 
relations represented by the two diagrams
\ref{compat} and \ref{compat2}.
The first case is just a calculation in the arc family picture. In 
the second case the arc family picture also
makes it very easy to write down
the family inducing the chain homotopy explicitly. 
We fix arc families
$a,b,c$ parameterized over $k_a,k_b,k_c$ respectively.

First notice that the family parameterized by $k_a  \times I 
\times k_b \times k_c$, which is depicted in figure \ref{compat}, illustrates
that
$$
a*(bc) \sim (a*b)c +(-1)^{|b|(|a|+1)} b(a*c)
$$
Second, the special two parameter family
shown in figure \ref{compat2}
--where the two parameter arc family realizing the chain homotopy is
parameterized over $k_a\times k_b \times T \times k_c$ with 
$T:=\{(s,t) \in [0,1] \times [0,1]: s+t \leq 1\}$--
gives the chain homotopy 
$$
(ab)*c \sim a(b*c) +(-1)^{|b|(|c|+1)}(a*c)b
$$
In both cases, we used normalizing homotopies as before.

These two equations imply that
\begin{equation}
\{a,b\cdot c\}\sim \{a,b\}\cdot c + (-1)^{|b|(|a|+1)}b\cdot\{a,b\},
\end{equation}
and we obtain:
\subsubth{Proposition} \sl The bracket $\{\; ,\;\}$ is
a Gerstenhaber bracket up to chain homotopy on $C_*(\Arc _{cp})$ or 
$C_*(\Arc)$ for the product $\cdot$.\rm\medskip

Summing up, we obtain:
\subth{Theorem} \sl There is a BV structure on $C_*(\Arc _{cp})$ or 
$C_*(\Arc)$ up to explicit chain boundaries.
The induced Gerstenhaber bracket is also given by such explicit boundaries.
This bracket is compatible with a product (associative and commutative
up to explicit boundaries) 
given by  a point in $C_0(\Arc_{cp}(2))$.\rm

\subth{Corollary} \sl All Arc$_{cp}$ and Arc algebras 
are Batalin-Vilkovisky algebras. The Gerstenhaber structure induced
 by the BV operator coincides with the bracket induced by the pre--Lie
product $*$. Hence all Arc$_{cp}$ and 
Arc algebras are GBV algebras.\rm

\section{Cacti as a suboperad of $\mathcal{ARC}$}
\label{cacti}

In the last section (\S \ref{BV}), we exhibited a GBV structure up to chain 
homotopy on $C_*(\Arc)$. Inspecting the arc families 
realizing the relevant homotopies, we observe
that the Gerstenhaber structure is already present in a 
suboperad which we shall discuss here.
This suboperad, called  ``linear trees'', corresponds
to the ``spineless cacti'' of \cite{K}. Furthermore for the BV operator, we 
only need to add the one more operation $\Delta$, so that the BV structure is
realized on the suboperad generated by spineless cacti and $\Arc_{cp}(1)$.
Lastly the suboperad generated by spineless cacti is contained
in the suboperad generated by cacti and the two Gerstenhaber structures
agree, i.e.\ the one coming from the BV operator and the previously defined 
bracket.

We will furthermore show that this operad is
the image of an embedding (up to homotopy) 
of Voronov cacti \cite{V} into $\Arc$.

By the results of \cite{CS,V}, algebras over the cacti operad have a 
BV structure. The map of operads we construct will
thus induce
the structure of BV algebras for algebras over the homology of $\Arc$.
This induced structure is indeed  
the same as that defined in section \ref{BV} as can
be seen from the embedding and our explicit realization of all
relevant operations.

\subsection{Suboperads}
\subsubth{Definition} The {\it trees suboperad} is defined for arc families
in surfaces with $g=s=0$ in the notation of \ref{rem4} by the allowed incidence matrix
$A^{(n)}$, whose non--zero
entries are $a_{0i}=1=a_{i0}$, for $i = 1, \dots, n$, and 
required incidence relations
$R^{(n)}$, whose entries are all equal to one.

This is a suboperad of $\Arc_{cp}$, and a representation of it as a collection of
labeled trees can be found in \cite{K}.

Dropping the requirement that $g=s=0$, we obtain a suboperad of $\Arc$
called the {\em rooted graphs} or {\em Chinese trees} suboperad.

\subsubth{Remark} We have already observed
that there is a linear and --- by forgetting the starting point
--- a cyclic order
on the set of arcs incident on each boundary component. In the 
(Chinese) trees suboperad all  bands must hit the $0$--th component, which induces
a linear and a cyclic order on all of the bands.
Furthermore the cyclic order is ``respected'' for trees,
 in the sense that the bands
meeting the $i$--th component form a cyclic subchain in the cyclic
order of all bands, for Chinese trees, this is an extra condition. 
And again all Chinese trees which satisfy this condition form a suboperad
which we call the {\em cyclic} Chinese trees.
The linear order is, however, not even respected for trees,
as can easily be seen in $\Arc_{cp}(1)$.

\subsubth{Linearity Condition} We say that an element
of the (cyclic Chinese) trees suboperad satisfies the
Linearity Condition if the linear orders match, i.e., the bands
hitting each boundary component in
their linear order are a subchain of all the bands in their linear
order derived from the $0$--th boundary.
 
It is easy to check that this condition is stable under composition.

We call the suboperad of elements
satisfying the Linearity Condition of the 
(cyclic Chinese) trees operad the {\em (cyclic Chinese) linear trees operad}.

\subsubth{Proposition} \sl The suboperad 
generated by (cyclic Chinese) linear
trees and $\Arc_{cp}(1)$ inside $\Arc$  coincides with
(cyclic Chinese) trees.\rm

\proof Given a (Chinese) tree
we can make it linear by gluing on twists from $\Arc(1)$ at the various boundary components 
as these twists have the effect of moving the marked point of the boundary
around the boundary. Since the cyclic order is already respected, such twists may be applied to arrange that the linear orders agree.  Since
(Chinese) linear 
trees and $\Arc(1)$ lie inside 
the (Chinese) trees operad the reverse inclusion is obvious.
\endproof

\subsection{Cacti}
There are several species of cacti, which are defined in 
\cite{K}, to which we refer the reader for details. 
By {\it cacti}, we mean  Voronov cacti as defined in \cite{V}, i.e.,
as  connected, planar tree-like configurations of 
parameterized loops, together with a marked point
on the configuration. This point, called ``global zero'', defines an
outside circle or perimeter by taking it to be the starting point and then
going around all loops in a counterclockwise fashion by jumping onto the
next loop (in the induced cyclic order) at the intersection points.
The {\it spineless} variety of cacti is obtained by postulating that 
the local zeros 
defined by the parameterizations of the loops coincide with the first
intersection point of the perimeter with a loop (sometimes called a ``lobe'')
 of the cactus. 

The gluing $\circ_i$ of two cacti
$C_1$ and $C_2$ is done by first scaling in such a way that the length of 
i--th lobe of $C_1$ matches the length of the outside circle of $C_2$
and then 
inserting the cactus into the i--th lobe by using the parameterization
as gluing data.

\subsubsection{Scaling of a cactus}
\label{cactscale} Cacti and spineless cacti both come with
a universal scaling operation of ${\mathbb R}_{>0}$
 which simultaneously 
scales all radii by the same factor $\l \in {\mathbb R}_{>0}$.
This action is a free action and the gluing descends to the
quotient by this action.

\subsubsection{Left, Right and Symmetric Cacti}
\label{leftrightsymmetric}
For the operadic gluings one has three basic possibilities
to scale in order to make the size of the outer loop of the 
cactus that is to be inserted match the size of the lobe into which
the insertion should be made.
\begin{enumerate}
\item Scale down the cactus which is to be inserted. This is the
original version -- we call it the right scaling version.
\item Scale up the cactus into  which will be inserted. 
We call it the left scaling version.
\item Scale both cacti. The one which is to be inserted by the size of 
the lobe into which it will be inserted and the cactus into which 
the insertion is going to be taking place by the size of the outer loop of the 
cactus which will be inserted. We call this it the symmetric scaling version.
\end{enumerate}

All of these versions are of course homotopy equivalent and in the quotient
operad of $cacti$ by overall scalings, the projective 
cacti $cacti/\mathbb{R}_{>0}$ they all descend to the same glueing.

\subsubsection{Framing of a spineless cactus}
We will give a map of spineless cacti into $\Arc$ called a framing.
 First notice that a spineless cactus can be decomposed by the initial point
and the intersection points into a sequence of arcs following the
natural orientation given by the data. These arcs are labeled by their
lengths as parts of parameterized unit circles.
To frame a given  spineless cactus, draw a pointed
circle around it and run an arc from each arc of the cactus
to the outside circle respecting the linear order, starting with the
initial arc of the cactus as the first arc emanating from the outside circle
in its orientation. Label each such arc by the parameter associated
to the arc of the cactus.

We can think of attaching wide bands to the arc of the cactus. 
The widths of the bands are just the lengths of the arcs to which they are attached.
Using these bands we identify the outside circle with the circumference
of the cactus. 
Notice that this ``outside'' circle appears in the
gluing formalism for cacti (see figures \ref{csops} to \ref{cactfig}). 

\subsubth{Proposition} \sl
\label{spinelessprop}
The operation of framing 
 gives a map of operads of spineless cacti  into $\Arc _{cp}\subset\Arc$
 whose image lies in the linear trees suboperad. 
Moreover we have a commutative square
$$
\begin{CD}
\text{spineless cacti}  @>frame>> \darc\\
@V\pi VV @VV\pi V\\
\text{spineless cacti}/{\mathbb R}_{>0}  @>\pi frame>> \Arc\\
\end{CD}
$$
where $\pi frame$ is defined by taking any lift of a 
spineless cactus/${\mathbb R}_{>0}$ to spineless cacti.
The map $frame$ becomes an operadic map, if one uses the
symmetric glueing for cacti. \rm 

\proof It is obvious that the image of the framing
of  spineless cacti lies inside the 
linear trees suboperad.
It is also clear that two cacti which differ by a scaling of
the type \ref{cactscale} get mapped to the same arc family, so that
the framing factors through spineless cacti modulo ${\mathbb R}_{>0}$.
The gluings  for cacti and arc
families are equivariant with respect to framing,
 as is obvious from the point of view of identifying
the outside circle with the circumference of a cactus.
By  considering weights rather than projective weights on the arcs of
a framing of a cactus, we can lift the framing of a cactus to $\darc$.
\endproof

To illustrate this map we will provide several 
examples.

\subsubsection{The basic operations of spineless cacti and their images}
\begin{figure}[ht!]
\epsfxsize = \textwidth
\epsfbox{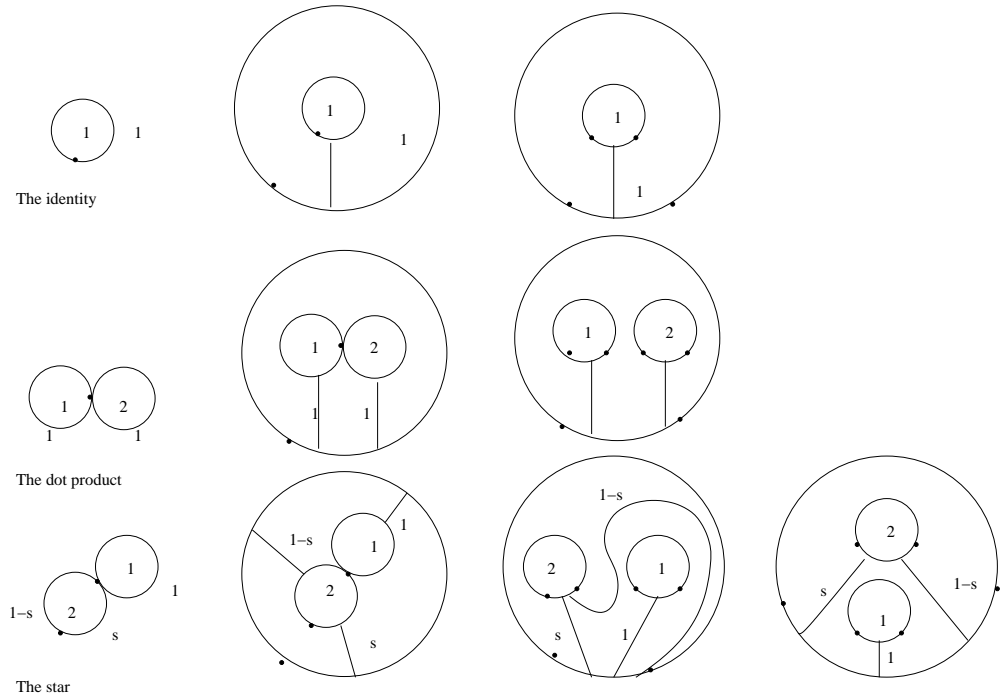}
\caption{\label{csops}The generators of spineless cacti up to normalization}
\end{figure}

 The operations: id, dot and $*$ are pictured in figure \ref{csops}
as well  as their images. This accounts for all 1 and 2 
boundary cacti up to an overall normalization.

\begin{figure}[ht!]
\epsfxsize = .85\textwidth
\cl{\epsfbox{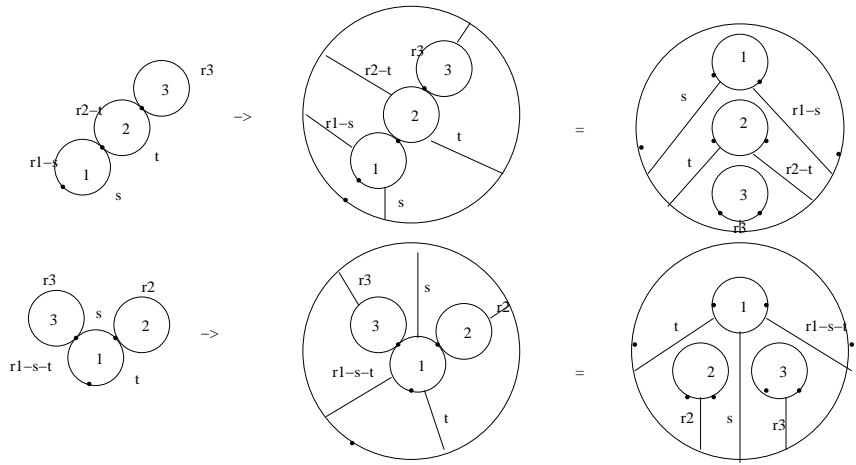}}
\caption{\label{threeops}The generic three boundary  spineless cacti}
\end{figure}

\begin{figure}[ht!]
\epsfxsize = .7\textwidth
\cl{\epsfbox{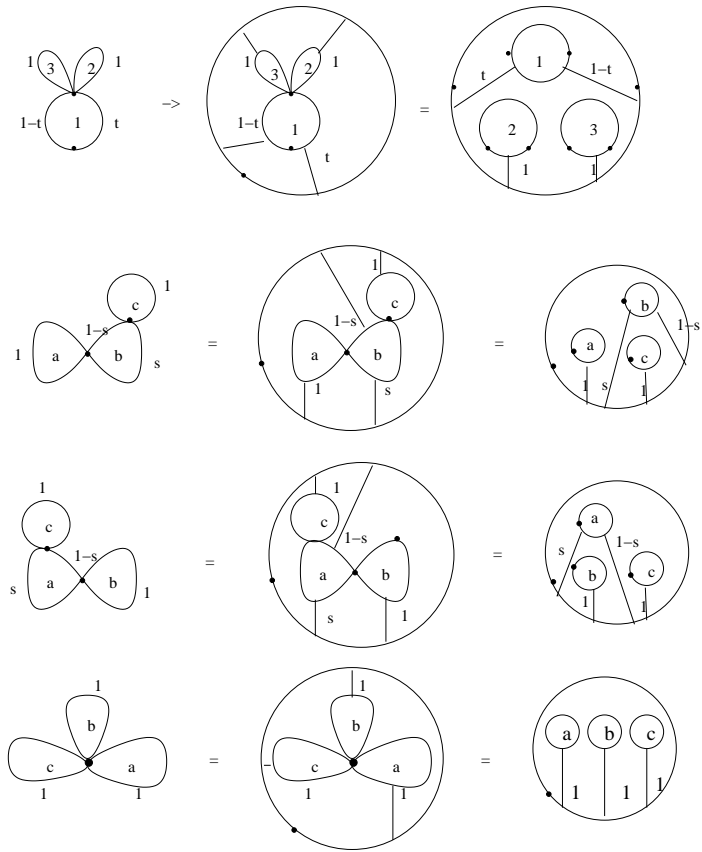}}
\caption{\label{threeops2} Special three boundary 
spineless cacti and their images}
\end{figure}

The generic cacti and their images are presented in figures \ref{threeops}
and \ref{threeops2}. All degenerate configurations are contained in the examples of figure 
\ref{threeops2}.

\begin{figure}[ht!]
\epsfxsize = \textwidth
\epsfbox{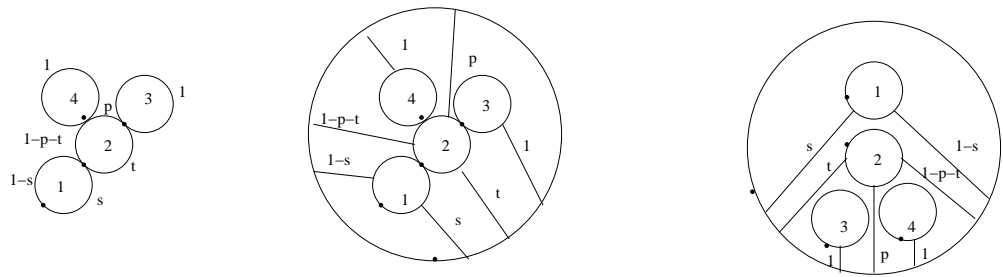}
 \caption{\label{cactfig} A  four boundary spineless cactus}
\end{figure}

\subsubth{Example} As a benefit, 
an example of a framing of a four boundary cactus is included in 
figure \ref{cactfig}.

\subsection{Framing of a cactus}
We wish to conclude this section with the framing procedure for cacti.
It is essentially the same as for spineless cacti. 
We decompose the cactus into the arcs of its
perimeter, where the break point of a cactus with spines are the
intersection points, the global zero and the local zeros.  One then
runs an arc from each arc to an outside pointed circle which is to
be drawn around the cactus configuration. The arcs should be
embedded starting in a counterclockwise fashion around the
perimeter of the circle. The marked points on the inside boundaries
correspond to the local zeros of the inside circles
viz.\ lobes of the cactus. 

\begin{figure}[ht!]
\epsfxsize = \textwidth
\epsfbox{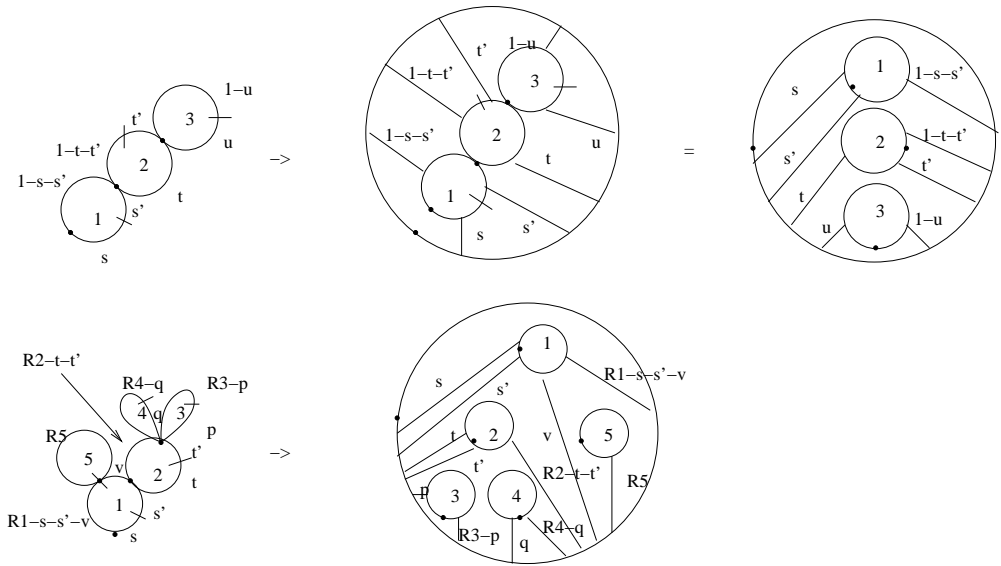}
\caption{\label{spinecacti}Framings of cacti}
\end{figure}

Two examples of this procedure are provided in figure \ref{spinecacti}.

Analogously to Proposition \ref{spinelessprop} we have:

\subsubth{Proposition}\sl
The operation of framing 
 gives a map of operads of cacti  into $\Arc _{cp}\subset\Arc$
whose image lies in the  trees suboperad. 
Moreover we have a commutative square
$$
\begin{CD}
\text{cacti}  @>frame>> \darc\\
@V\pi VV @VV\pi V\\
\text{cacti}/{\mathbb R}_{>0}  @>\pi frame>> \Arc\\
\end{CD}
$$
where $\pi frame$ is defined by taking any lift of a 
 cactus/${\mathbb R}_{>0}$ to cacti and we again regard
symmetric glueing for cacti.\rm

\subth{Remark} In the previous section, the arc families we employed
to provide the homotopies to exhibit the GBV structure had another 
common characteristic: their total 
weights at all boundary components except at the $0$--th boundary component
were all equal (to one). This property is not stable under gluing though.
However, one can define a quasi--operad of normalized cacti \cite{K}, where a
quasi--operad is an operad which is not necessarily associative.
Spineless cacti are however homotopy associative as an operad and
moreover there is a chain decomposition for spineless cacti
whose cellular chains give
an honest operad under the induced compositions \cite{K}. Furthermore this 
cellular operad gives a solution to Deligne's conjecture on the Hochschild
cohomology of an associative algebra \cite{K2}.
Under
framing normalized spineless cacti are a 
quasi--suboperad up to homotopy, i.e., the glueings are not
equivariant, but are equivariant up to homotopy. 
In fact in \cite{K} it is shown that 
(spineless) cacti are the direct product of normalized (spineless) cacti 
and a scaling operad built on ${\mathbb R}_{>0}$, and moreover
all structures such as the Gerstenhaber and BV can be defined on the level
of normalized cacti.

\section{The Loop of an arc family}
\label{loop}

Given a surface with arcs we can forget some of the structure and
in this way either produce a collection of loops or one loop which
is given by using the arcs as an equivalence relation. Using these
two maps we obtain an operation of the Chinese trees suboperad of
 $\Arc$ on the loop space of any manifold on the homological level.

We expect that this action can be enlarged to all of $\Arc$,
but for this we would need a generalization of the results of 
Cohen and Jones \cite{CJ}, which we intend to study elsewhere.

\subsubsection{The boundary circles}
Given an exhaustive weighted arc family $(\a )$ in the surface $F$, we can consider the measure-preserving
maps
\begin{equation}
\tilde {c}^{(\a)}_i :\del_i(\a) \rightarrow S^1_{m_i}
\end{equation}
where $S^1_r$ is a circle of radius $r$ and $m_i= \mu^i(\del_i(\a))$ is the total
weight of the arc family  at the $i$--th boundary. 
Combining these maps, we obtain
\begin{equation}
\tilde c: \del(\a)\rightarrow \coprod_i S^1_{m_i}.
\end{equation}
Choosing a measure on $\partial F$ as in \S 1 to identify $\partial (\a)$ with
$\partial F$, we finally
obtain a map 
\begin{equation}
circ: \del F \rightarrow \coprod_i S^1_{m_i}
\end{equation} 
Notice that the image of the initial points of the bands give well-defined
basepoints $0 \in
S^1_{m_i}$ for each
$i$.

\subsection{The equivalence relations induced by arcs}

On the set $\del(\a)$ there is a natural reflexive and symmetric relation
given by  $p \sim_{fol}q$ if $p$ and $q$ are on the same leaf of the partial measured foliation.

\subth{Definition}
Let $\sim$ be the equivalence relation on $\coprod_i S^1_{m_i}$
generated by $\sim_{fol}$.
In other words $p \sim q$ if there are leaves $l_j$, for $j= 1,\dots m$, so that
$p \in \tilde c(\del(l_1)), q\in \tilde c(\del(l_n)) $ and 
$\tilde c(\del(l_j)) \cap \tilde c(\del(l_j+1))\neq\emptyset$.

\subsubth{Remark} It is clear that neither the image of $circ$
--- which will denote by $circ((\a))$ --- as
a collection of parameterized
circles nor the relation $\sim$ depends upon the choice of measure
on $\partial F$.

\subth{Definition}
Given a deprojectivized arc 
family $(\a) \in \darc$, 
we define $Loop((\a))=circ((\a))/\sim$ and denote the projection map
$\pi:circ((\a))
\rightarrow Loop((\a)) $.

Furthermore, we define two maps taking
values in the monoidal category of pointed spaces:
\begin{eqnarray}
int((\a)) &=& \bigsqcup_{i=1}^n (\pi(\tilde c_{i}^{(\a)}(\del_i (\a)),\pi(*_i))\\
ext((\a)) &=& ( (\pi(\tilde c_{0}^{(\a)}(\del_0 (\a)),\pi(*_0))
\end{eqnarray}
and call them the internal and external loops of $(\a)$ in $Loop((\a))$.
We denote the space with induced topology given by the collection
of images $Loop((\a))$ of all $(\a) \in \darc(n)$ by $Loop(n)$.

Notice that there are $n+1$ marked points on $Loop((\a))$
for $(\a) \in \darc(n)$

\begin{figure}[ht!]
\epsfxsize = 4in
\cl{\epsfbox{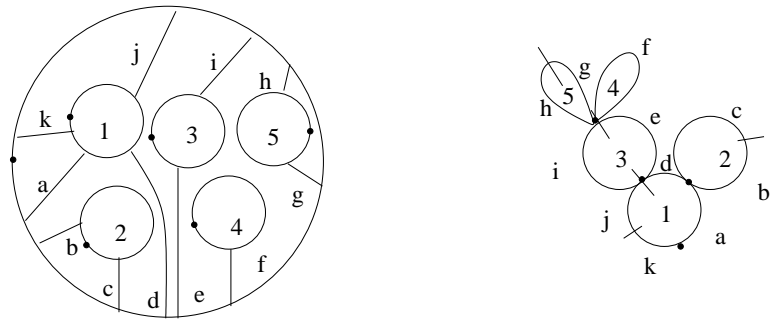}}
\caption{\label{loop1}
An arc family whose loop is a cactus}
\end{figure}

\begin{figure}[ht!]
\epsfxsize = 3in
\cl{\epsfbox{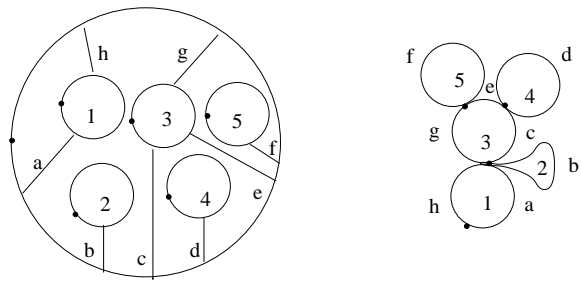}}
\caption{\label{loop2} An arc family whose loop is a cactus without spines}
\end{figure}

\begin{figure}[ht!]
\epsfxsize = 4.2in
\cl{\epsfbox{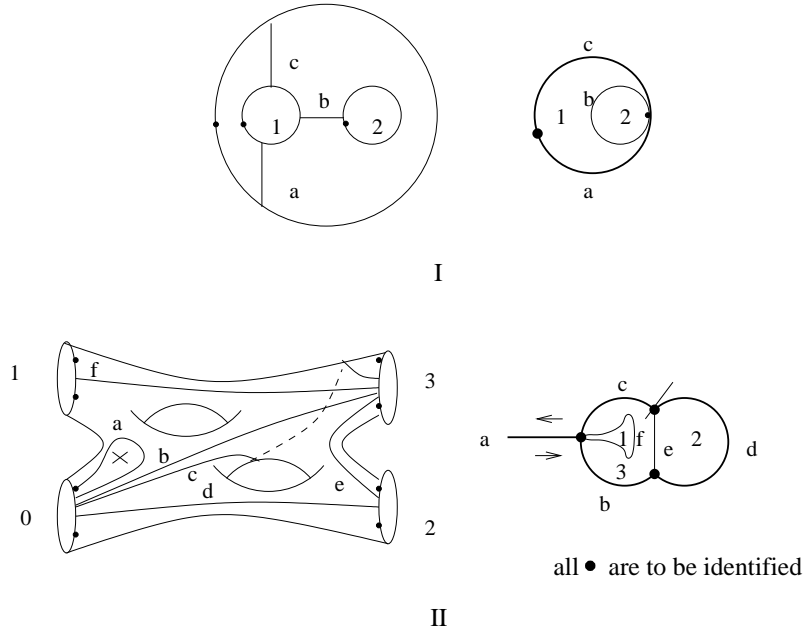}}

\caption{\label{loop3}
Loops of arc families not yielding cacti: I, genus 0 case,\qua
II, genus 2 with one puncture}
\end{figure}

Examples of loops of an arc family are depicted in figures
\ref{loop1}--\ref{loop3}. In figure \ref{loop3} I the image of the
boundary 1 runs along the outside circle and then around the inside
circle. The same holds for the boundary 3 in figure \ref{loop3} II.
In both \ref{loop3} I and II, the outside circle and 
its basepoint are in bold.


\subsubth{Remark} There are two types of intersection points
for pairs of loops. 
The first are those coming from the interiors
of the bands; these points are double points and occur 
along entire intervals.
The second type of multiple point arises from the boundaries of the
bands via the transitive closure; they can have any multiplicity
but are isolated.

\subsection{From loops to arcs}
If the underlying surface of an arc family satisfies $g=s=0$, then its Loop together
with the parameterizations uniquely determines the arc family.
In other words, the map $frame$ is a section of $Loop$.

\subsubth{Definition} A {\it configuration of circles} is the image of a surjection
$p:\coprod_i S^1_{m_i} \rightarrow L$ of metric spaces such that each point of $L$ lies in the
image of at least two components and the intersections of the images
of more than two components are isolated.
Let $Config(n)$ be the space of all such configurations
of $n+1$ circles with the natural topology.
We call a configuration of circles {\em planar}, if $L$ can be
embedded in the plane with the natural orientation for all images 
$S^1_i:i\neq 0$ coinciding with
the induced orientation and the opposite orientation for $S^1_0$.
We call the space of planar configurations of $n+1$ circles
$Config_{p}(n)$.

\subsubth{Proposition}\sl
The map $Loop:\darc(n)\rightarrow Config_{p}(n)$ is surjective.\rm

\proof To describe a right inverse, fix 
$(p:=\coprod_i S^1_{m_i} \rightarrow L) \in Config_{p}(n)$,
fix an embedding of $L$ and 
$F_{0,n+1}^0$ in the plane, 
fix an identification of each boundary component with $S^1$, and identify
the disjoint union of these boundary components regarded as labeled 
circles with the source
of $p$. 
On $L$, remove the images of the basepoints. For 
each component of dimension one in the intersection of two of the 
components of the map $p$, draw an arc between the respective boundaries 
on $F_{0,n+1}^{0}$ of width given by the length of that component.
These arcs are to be embedded in the linear order dictated by the various 
parameterizations and incidence conditions in the plane. If it happens
in this way that an arc arises that is parallel to the boundary,
then we insert a puncture in the region between the arc and the respective
boundary.  
Since the higher order intersections are isolated
they do not contribute. It is easily checked that this yields a right inverse of Loop.
\endproof

\subsubth{Proposition and Definition} \sl The deprojectivized 
arc families such that  
$\pi|_{\tilde c(\del_0)}(\a))=Loop((\a))$  
constitute a suboperad of $\Arc$. We call this suboperad {\em $\Loop$}.\rm 

\subsubth{Proposition}\sl  If $(\a)\in \Loop$ then $Loop((\a))$ is a
cactus.  Furthermore, the operad $\Loop$ is identical to the operad
of Chinese trees.\rm

\proof If  $\tilde c(\del_0)(\a)\supset Loop((\a))$, then 
there are no arcs running between two boundary components if neither is $\partial _0$.
Furthermore, since $\pi|_{\del_0}(F)\subset \tilde c(\del_0)(\a)$, there is no
arc running from $\del_0$ to $\del_0$.\endproof
Collecting the results above, we have shown:

\subth{Theorem} \sl The framing of a  cactus
 is a section of $Loop$ and 
is thus an  embedding. This embedding identifies
(normalized and/or spi\-ne\-less) cacti as (normalized and/or linear) trees.
$$
\begin{CD}
@<Loop<<\\
\text{(spineless) cacti}  @>>frame> \darc\\
@V\pi VV @VV\pi V\\
\text{(spineless) cacti}/{\mathbb R}_{>0}  @>\pi frame>> \Arc\\
@<<\pi Loop<
\end{CD}
$$
where $\pi Loop$ is defined by choosing any lift and the glueing
for cacti is the symmetric glueing.
\rm

\subsection{Comments on an action on the loop space}
Given a manifold $M$ we can consider its loop space $LM$. Using
the configuration we have maps
\begin{equation}
Arc(n) \times LM^n\stackrel{Loop \times id}
{\longleftarrow} Config(n)\times LM^n
\stackrel{i}{\longleftarrow}{L}^{Config(n)}M
\stackrel{e}{\rightarrow} LM
\end{equation}
where ${L}^{Config(n)}M$ are continuous 
maps of the images $L$ of the configurations
into $M$, i.e., such a map takes a configuration 
$p:\coprod_{i} S^1_{m_i} \rightarrow L$
and produces a continuous $f:L \rightarrow M$; 
the maps $i,e$ are given by
$i(f)=(p:\coprod_{i} S^1_{m_i} \rightarrow L,f(p(S^1_{m_1}),\dots,
f(p(S^1_{m_n})))$ and 
$e(f)=f(p(S^1_{m_0}))$.

One would like to follow the Pontrjagin--Thom construction of \cite{CJ,V} so that
the maps $i,e$ in turn would induce maps on the level of homology
\begin{multline}
``H_*(\Arc(n)) \otimes H_*(LM^n)\simeq
H_*(\darc(n)) \otimes H_*(LM^n)
\stackrel{Loop_*}{\rightarrow}\\
H_*(Config(n)) \otimes H_*(LM^n)
\stackrel{i^!}{\rightarrow}H_*({L}^{Config(n)}M)
\stackrel{e_*}{\rightarrow} H_*(LM)\text{''}
\end{multline}
where $i^!$ is the ``Umkehr'' map, {\it but the
map $i^!$ is only well defined on the subspace of cacti}.

If we restrict ourselves to this subspace we obtain from 
the above following \cite{CS,V,CJ}:

\subsubth{Proposition}\sl  
The homology of the loop space of a manifold 
is an algebra  over the suboperad of cyclic Chinese trees.\rm 

\subsubsection{Remarks} 
\begin{figure}[ht!]
\epsfxsize = \textwidth
\epsfbox{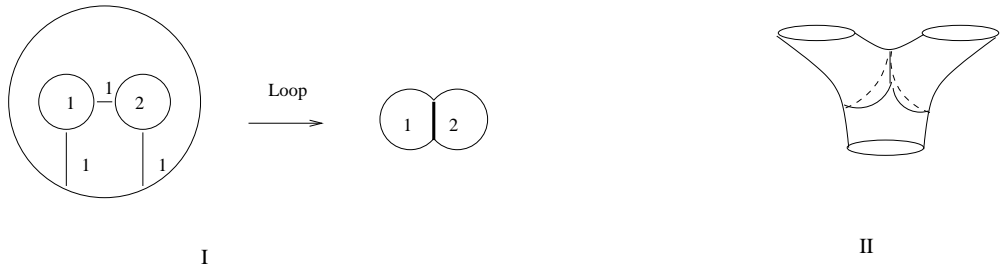}
\caption{\label{zwie}
I, the Loop of a symmetric pair of pants,\qua 
II, a closed string field theory picture of a pair of pants}
\end{figure}

\begin{itemize}
\item[(1)]
It is clear that one desideratum is the
extension of this result to all of $\Arc$. 
\item[(2)]
The first example of an operation of composing loops which are
not cacti would be given by the $Loop$ of the pair of pants with
three arcs as depicted in figure \ref{zwie}. This kind of composition
first appeared in the considerations of closed string field theory.
\item[(3)] If the image of Loop is not connected, then the information is 
partially lost. This can be refined however by using a prop version 
of our operad.
\item[(4)] Factoring the operation through Loop has the effect that 
the internal topological structure is forgotten;
thus, the  torus with two boundary components has the same effect as the cylinder
for instance.
\end{itemize}

\section{Twisted Arc Operads}
Recall from Remark \ref{remnu} that given an arc family $[\a]$, we may identify each of the  ends of
the bands with $S^{1}$, or if given some additional data also with a boundary component of $F$. We will use this geometric fact in
this section to geometrically combine the circle operads defined in the
Appendix with the Arc operad, in the sense that we can identify
the boundary of a surface with $S^1$ via the maps \ref{circleboundary}.
This leads to several direct and semi--direct (in the sense of \cite{K}) 
products.  

\subth{Definition and Proposition} We define the differential 
$Arc$ operad
$dArc$ to be the Cartesian product of $\Arc$ with the operad $d$
where we think of the $S^{1}$s of $d$ as being the boundaries
$1,\ldots,n$.
Algebras over this operad are $Arc$ algebras together with a
differential $d$ of degree one that is also a derivation and that
anti--commutes with $\delta$.

\subsubth{Remark} By the previous sections and the Appendix, 
we see that algebras
over this operad will have the structure of dGBV algebras. This structure
will already be present for algebras over the suboperad of (Chinese) trees
and algebras over the suboperad of (linear and/or Chinese) trees
will have the structure of differential Gerstenhaber algebras.

\subth{Definition and Proposition} We define the 
untwisted ${\lr}Arc$ operad
to be the Cartesian product of $\Arc$ with the operad $\lr$
where we think of the $S^{1}$s of $\lr$ as being the boundaries
$0,\ldots,n$. This operad is cyclic.
Algebras over this operad are Arc algebras together with two more
differentials $\ldel$ and $\rdel$ of degree one that satisfy the
relations of Theorem \ref{grcyclicsphere}. Both $\ldel$ and $\rdel$ anti--commute with
$\delta$.

\subth{Remark} There are no semi--direct products with these operads,
at least none that are strictly associative. One can produce, however,
by using the standard procedures \cite{SW,MSS} and the gluing below operads
that are associative up to homotopy.

\subth{Definition} We set $\widetilde{Arc}_{g,s}(n):=
Arc_{g,s}(n) \times (S^{1})^{n+1}$.
Recall that the gluing \ref{glue} was achieved by identifying the
boundaries of the bands with the help of the identification of each
of the two boundaries with $S^1$ and then identifying these $S^1$s.
This construction can be twisted by using the natural $S^1$ action on
one or both of the $S^1$s before making the identification.

We 
use this freedom to define twisted gluings depending on one
variable $\th \in S^1$ (recall that we use additive notation):  Define
$\circ_i^{\th}$ to be the gluing resulting from the identifications
$c^{\a}_i: \del_i(\a) \rightarrow S^1$ and $s_\th \circ
c^{\b}_0:\del_0(\b)\rightarrow S^1$, where $s_th (\phi) =\phi +\th$.

Now we introduce operations
$$\circ_{i}: \widetilde{Arc}_{g,s}(n) \times \widetilde{Arc}_{g',r}(m)
\rightarrow \widetilde {Arc}_{g+g',r+s}(m+n-1)$$ 
as follows: Let
$\vec{\th}=(\th_0,\dots,\th_n)$ and $\vec{\th'}=(\th'_0,\dots,\th'_m)$
\begin{equation}
([\a'],\vec{\th})\circ_i([\b'],\vec{\th'}):= (\a \circ_i^{-\th'_0}
\b, \vec{\th}\circ_i \vec{\th'})
\end{equation}
where $ \vec{\th}\circ_i \vec{\th'}$ are the operations of $\rr$.

The geometric meaning of this is that we regard the points in $S^1$
as points on the boundary of the band and use the point on the boundary component 0
as a possible offset for gluing. Upon gluing we use the point on the
boundary component $0$ of the family $\b$ to be the beginning of the
window instead of the point $0$. To implement this on the whole
surface, we use the diagonal shift action. Lastly, we use the offset
of the i--th boundary of $\a$ to translate all the offsets of $\b$.

\subsubth{Theorem} \sl The operations $\circ_{i}$ turn the
$\mathbb S$--module $\widetilde{Arc}_{g,s}$ into an operad.\rm

\proof Again the operad structures are evident. This is due to the
fact that the twisted gluing of the bands is dependent only on $\theta_0$
which does not play a role in the $\rr$ gluing.\endproof

\subth{Remark} In \cite{K} we also defined bi--crossed products
of operads, which provide the right framework for the cacti without
spines operad relative to the cacti operad. This structure can be carried over
to $\Arc$ by looking at the action of $S^1$ moving the marked point 
of the boundary zero. Taking the cue from \cite{K}, we would say that for
each boundary connected to the $0$--th boundary
there is a natural parameterization induced by marking the first point of the
first band in the window of the $i$-th boundary if it is 
connected to the $0$-th boundary. This identification will induce
a bi--crossed product of the type \cite{K}, which now only acts partially
and not through the whole diagonal.

\renewcommand{\theequation}{A-\arabic{equation}}
\renewcommand{\thesection}{A}
% redefine the command that creates the equation no.
\setcounter{equation}{0}  % reset counter 
\setcounter{subsection}{0}
\section*{Appendix:\newline Circle operads, differentials and derivations}  
% use *-form to suppress numbering

In this appendix, we study operads concocted out of Cartesian
products of circles with basepoint whose algebras over the respective
homology operads will be graded commutative and associative algebras
together with several different types of operators,
which we used to build several direct and semi--direct products
of these operads with the $\Arc$ operad.
In the process we 
give a classification of all linear and 
local operads built on $(S^1)^{n+1}$.

We will view the circle $S^1$ as ${\mathbb R}/{\mathbb Z}$, so that it has
basepoint
$0$ and has a natural additive angular coordinate $\theta$.

\subsection{The Operad $(S^1)^{n}$} Although this operad
seems to be known to the experts (and appears in the semi--direct products
of operads with groups \cite{SW,MSS}), we shall
give its complete description here. One interesting point is that
its homology operad is that of the operad built on the supervector space
$\ZZ$ as defined in \cite{K}.

\subsubth{Definition} Let $d$ be the collection of spaces
$d(n) := (S^1)^{ n}$ together with the operations
$\circ_i: d(n) \times d(m) \rightarrow d(n+m-1)$ defined by
$$(\theta_1,\dots,\th_n)\circ_i(\th'_1, \ldots, \th'_m)=
(\th_1, \ldots, \th_{i-1}, \th'_1+ \th_i, \ldots, \th'_m +\th_i, \th_{i+1},
\dots, \th_n)$$ and  the natural permutation action of $\Sn$.

\subsubth{Proposition}\sl $d$ is a topological operad.
The algebras over the homology operad $H_{*}(d)$ are differential
graded, associative, commutative algebras. \rm

\proof The axioms for an operad are straight-forward to verify.
The unit is the class of $0$ in $S^1$. By the K\"unneth formula, we have
$H_{p}(d(n))= {\mathbb Z}^{n\choose p}$. Hence,
$H_*(d(n))$ is a free ${\mathbb Z}$-graded module, and a basis can be indexed
by
$(\Zz)^{n}$. Using this basis, the $\circ_{k}$'s read:
\begin{align*}
&(\alpha_{1},\ldots,\alpha_{n}) \circ_{k}(\beta_{1},\ldots,\beta_{m}) =\\
&\pm (1-\alpha_k)(\alpha_{1},\ldots,\alpha_{k-1},\beta_{1},\ldots, \beta_{m},
\alpha_{k+1},\ldots,\alpha_{n}) + \\
&\alpha_k \sum_{l=1}^m \pm (\beta_l+\alpha_k)[2]
(\alpha_{1},\ldots,\alpha_{k-1},\beta_{1},\ldots,\beta_l+\alpha_k,\ldots, \beta_{m},
\alpha_{k+1},\ldots,\alpha_{n})
\end{align*}
where $\pm$ is the standard supersign, and $[2]$ means modulo 2.

Consider the operad ${\cal P}$ generated by the $\SS$-module $E=\del\SS_1\oplus
\mu\cdot k$, where $k$ is the trivial $\SS_2$
representation, together with the relations $\del\del:=\del\circ_1\del=0$,
$\del\mu=\mu\circ_1\del+\mu\circ_2\del$ and $\mu\circ_1\mu=\mu\circ_2\mu$.
We easily check that there is an isomorphism
$$\psi : {\cal P}\rightarrow H_*(d)$$
defined by $\psi(\del)=(1)$ and $\psi(\mu)=(0,0)$. The inverse
morphism $$\gamma : H_*(d)\rightarrow {\cal P}$$
is given by $\gamma(\alpha_1,\ldots,\alpha_n)=\mu^{(n)}(\del^{\alpha_1},\ldots,
\del^{\alpha_n})$, where $\mu^{(n)}$ is the $(n-1)$-fold composition of
$\mu$ for $n\geq 1$.

As a consequence, an algebra over $H_*(d)$ is an algebra over ${\cal P}$,
that is, a graded vector space ${\cal A}$
together with an operation $\del$ of degree $1$ and a graded
commutative multiplication $\cdot$
of degree $0$ satisfying $\del^2=0$, $\del(a\cdot b)=\del(a)\cdot b +
(-1)^{|a|} a\cdot \del(b)$, where $\cdot$ is associative.\endproof

\subsection{Classification of the operad structures on  $(S^1)^{ n+1}$}
Let $\cal P(n)=(S^1)^{ n+1}$.
This vector space is endowed with an action of the symmetric group
$\SS_n$ on its last $n$ variables.
The aim of this section is to give a complete classification of the
different kinds of operad compositions in $\cal P$,
where we demand only that the composition be local  and linear. This
classification
is given in the next proposition, whose straight-forward proof is
completed by explicating the various associativity
conditions.

\subsubth{Proposition}\sl \label{classificationS1} 
Let $\alpha,\beta,\gamma,\delta \in \R$
and define the composition on  $(S^1)^{n+1}$ by
$$\circ_i: (S^{1})^{n+1}\times (S^{1})^{m+1}
\rightarrow (S^{1})^{n+m}$$
\begin{multline*}
(\theta_0,\dots,\th_n)\circ_i(\th'_0, \ldots, \th'_m) :=\\
(\th_0 +\gamma \th_i+\delta \th'_{0}, \ldots, \th_{i-1}+\gamma
\th_i+\delta \th'_{0},\\
\th'_1+ \alpha \th_i+\beta \th'_0, \ldots, \th'_m +\alpha \th_i+\beta \th'_0,\\
\th_{i+1}+\gamma \th_i+\delta \th'_{0},\dots, \th_n+\gamma\th_i+\delta
\th'_{0})
\end{multline*}
There is a translation action of $\Z$ on the parameters,
which leaves the operations invariant, so 
the parameters may be taken in $\R/\Z$.

Then this composition endows $(S^{1})^{ n+1}$ with a structure
of operad only in the following cases:

\begin{itemize}
\item[] $\gamma=0$\qua and
\item[\rm(i)] $\alpha=\beta=\delta=0\vrule width0pt height 15pt$, or
\item[\rm(ii)] $\alpha=\delta=1$ and $\beta=0$, or
\item[\rm(iii)] $\alpha=1$ and $\delta=0$.
\end{itemize}
\rm

\subsubth{Remark} The first case, which corresponds to the cyclic 
operad
built on $S^1$ as a space (cf.\ \cite{K}),
will be described in next
proposition. The second case gives a cyclic and unitary operad, which
will be denoted by $\lr$ and is studied in next sub-section. The third
case is neither cyclic nor unitary but is of interest in combination with
the Arc operad; it will be denoted by $\rr_\lambda$, $\lambda=\beta\in\R$.

Corresponding to case i) of Proposition \ref{classificationS1}, we have:

\subsubth{Proposition} \sl Let $\cal Q$ be the operad defined by
\hfill\break
$\cal Q(n)=(S^1)^{n+1}$ together with the composition

$(\th_0,\ldots,\th_n)\circ_i(\th'_0,\ldots\th'_m)=(\th_0,\ldots\th_{i-1},
\th'_1,\ldots,\th'_m,\th_{i+1},\ldots\th_n)$.

\noindent Then an algebra over the homology operad $H_*(\cal Q)$ is a
graded commutative and associative algebra $A$ together with an operator
$\Delta : A\rightarrow A$ of degree 0, and two differentials
$\ldel,\rdel : A\rightarrow A$ of degree 1,
satisfying the following relations:
\begin{eqnarray*}
\Delta^2=\Delta,& \Delta\rdel=\rdel,& \ldel\Delta=\ldel,\\
\Delta\ldel=0,& \rdel\Delta=0,& \rdel\ldel=0,\\
\Delta(ab)=&\Delta(a)b=a\Delta(b)=&ab.\\
\end{eqnarray*}
\rm
\proof The proof (which is left to the reader) consists of
describing the operad
$H_*(\cal Q)$ by generators and relations as in the proof of Proposition
\ref{homcyclicsphere}. \endproof

\subsection{The Cyclic Operad $(S^1)^{n+1}$}
As seen in Proposition \ref{classificationS1}, there is only one
cyclic operad built on $(S^1)^{n+1}$ which is moreover
unitary. This cyclic operad is denoted by $\lr$, and the unit in $\lr(1)$
is given by 
$(0,0)\in (S^{1})^{ 2}$.

Recall that the composition is given by

$\circ_i: \lr(n) \times \lr(m) \rightarrow \lr(n+m-1)$
\begin{multline*}
(\theta_0,\dots,\th_n)\circ_i(\th'_0, \ldots, \th'_m) :=\\
(\th_0 +\th'_{0}, \ldots, \th_{i-1}+\th'_{0},
\th'_1+ \th_i, \ldots, \th'_m +\th_i,
\th_{i+1}+\th'_{0},\dots, \th_n+\th'_{0}),
\end{multline*}

\noindent so cyclicity is clear from the symmetric nature of the
operation. Nevertheless, we will write out the cyclicity calculation  once
explicitly. Denoting by $^{*}$ the action of the long cycle $(0 \cdots
n)\in \Snn$, we have

\begin{multline*}
[(\theta_0,\dots,\th_n)\circ_n(\th'_0, \ldots, \th'_m)]^{*} =\\
(\th_0 +\th'_{0}, \ldots, \th_{n-1}+\th'_{0},
\th'_1+ \th_n, \ldots, \th'_m +\th_n)^{*}=\\
(\th'_m +\th_n,\th_0 +\th'_{0}, \ldots, \th_{n-1}+\th'_{0},
\th'_1+ \th_n, \ldots, \th'_{m-1} +\th_n)
\end{multline*}
and
\begin{multline*}
(\th'_0, \ldots, \th'_m)^{*}\circ_1(\theta_0,\dots,\th_n)^{*} =
(\th'_m,\th'_0, \ldots, \th'_{m-1})\circ_1(\th_{n},\theta_0,\dots,\th_{n-1}) \\
=(\th'_m +\th_n,\th_0 +\th'_{0}, \ldots, \th_{n-1}+\th'_{0},
\th'_1+ \th_n, \ldots, \th'_{m-1} +\th_n)
\end{multline*}
as required.\endproof

The homology analogue is given in
Proposition \ref{homcyclicsphere}. Theorem
\ref{grcyclicsphere} gives the description of the operad $H_*(\lr)$ in
terms of generators and relations, whereas Corollary
\ref{corcyclicsphere} gives the description of algebras over the operad
$H_*(\lr)$.

\subsubth{Proposition}\sl \label{homcyclicsphere}
The homology of $\lr$ is a cyclic
operad.
It  is given by
$H_{p}(\lr(n))= {\mathbb Z}^{n+1\choose p}$. A basis of the free ${\mathbb Z}$-graded module
$H_p(\lr(n))$ is given by a sequence
$(\alpha_0,\ldots,\alpha_n)\in (\Zz)^{n+1}$, with
$\sum_{i=0}^{n}\alpha_i=p$. Denote by $e_k, 0\leq k\leq n+m-1$ the
canonical basis of $(\Zz)^{n+m}$. Then the composition
$$\circ_k : H_p(\lr(n)) \otimes H_q(\lr(m)) \rightarrow
H_{p+q}(\lr(n+m-1))$$
is given by
\begin{align*}
&(\alpha_0,\ldots,\alpha_n)\circ_k(\beta_0,\ldots,\beta_m) =
\pm (1-\alpha_k)(1-\beta_0) \Phi +\\
&(1-\alpha_k)\beta_0 \sum_{l=0}^{k-1} \pm (\alpha_l+\beta_0) [2]
(\Phi + \beta_0 e_l) + \\
&(1-\alpha_k)\beta_0\sum_{l=k+1}^{n} \pm (\alpha_l+\beta_0)[2]
(\Phi + \beta_0 e_{l+m-1})+\\
&\alpha_k(1-\beta_0) \sum_{p=1}^m \pm (\alpha_k+\beta_p) [2]
(\Phi+\alpha_k e_{k+p-1})+\\
&\alpha_k\beta_0 \sum_{p=1}^m \sum_{l=0}^{k-1}\pm
(\alpha_l+\beta_0)(\alpha_k+\beta_p)[2]
(\Phi+\beta_0 e_l+\alpha_k e_{k+p-1})+ \\
&\alpha_k\beta_0 \sum_{p=1}^m \sum_{l=k+1}^{n} \pm
(\alpha_l+\beta_0)(\alpha_k+\beta_p)[2]
(\Phi+\beta_0 e_{l+m-1}+\alpha_k e_{k+p-1})
\end{align*}
where $\Phi:=(\alpha_0,\ldots,\alpha_{k-1},
\beta_1,\ldots,\beta_m,\alpha_{k+1},\ldots,\alpha_n)$.
\rm

The next result gives a description of the homology
operad
$H_*(\lr(n))$ in terms of generators and relations.

\subsubth{Theorem}\sl \label{grcyclicsphere} Let $E_*$ be the following
graded
$\SS-$module:  $E_1$ is generated by
$\ldel$ and $\rdel$ in degree $1$,   $E_2$
is generated by $\mu$ in degree $0$, where the action of the cycle
$(12)\in\SS_2$ is given by $(12)\cdot \mu=\mu$ and $E_n=0$ for $n>2$.
Let ${\cal F}(E)$ be the free graded operad generated by $E$. Let
${\cal R}$ be the sub $\SS$-module of ${\cal F}(E)$ generated by the
elements
\begin{align*}
&\ldel\circ_1\ldel,\ \rdel\circ_1\rdel,\\
&\ldel\circ_1\rdel+\rdel\circ_1\ldel ,\\
&\rdel\circ_1\mu-\mu\circ_1\rdel-\mu\circ_2\rdel ,\\
&\mu\circ_1\ldel-\ldel\circ_1\mu-\mu\circ_2\rdel, \\
&\mu\circ_1\mu-\mu\circ_2\mu .
\end{align*}
Then there is an isomorphism between the operads
$${\cal F}(E)/<{\cal R}>\  \simeq \ H_*(\lr)$$
\rm

\proof Let $\psi : E\rightarrow H_*(\lr)$ be the $\SS$-module morphism
given
by $\psi(\ldel)=(1,0),\ \psi(\rdel)=(0,1)$ and $\psi(\mu)=(0,0,0)$,
which induces a morphism $\psi : {\cal F}(E)\rightarrow H_*(\lr)$. To
define a morphism $$\psi : {\cal F}(E)/<{\cal R}> \rightarrow H_*(\lr),$$
we have to prove that $\psi({\cal R})=0$. This is a straight-forward
computation; as an example, we compute
\begin{align*}
&\psi(\mu\circ_1\ldel-\ldel\circ_1\mu-\mu\circ_2\rdel)\\
&=(0,0,0)\circ_1(1,0)-
(1,0)\circ_1(0,0,0)-(0,0,0)\circ_2(0,1)\\
&=(1,0,0)+(0,0,1)-(1,0,0)-(0,0,1)=0.
\end{align*}
The next step of the proof is to define an inverse 
$$\gamma : \lr\rightarrow {\cal F}(E)/<{\cal R}>$$ to $\phi$.
Denote by $\mu^{(n)}$ or $\mu$ the $(n-1)$-fold composition of
$\mu$ which is independent of the manner of composition since $\mu$ is
associative, by definition of ${\cal R}$.
The following statements are immediate by induction on $n$
\begin{equation}\label{equardel}
\rdel\circ_1\mu=\sum_{i=1}^n \mu\circ_i \rdel
\end{equation}
\begin{equation}\label{equaldel}
\mu\circ_i\ldel=\ldel\circ_1\mu+\sum_{j\not=i} \mu\circ_j\rdel
\end{equation}
Let $\gamma : H_*(\lr) \rightarrow {\cal F}(E)/<{\cal R}>$ be the map
defined by
$$\gamma(\alpha_0,\ldots,\alpha_n)=
\ldel^{\alpha_0}\mu^{(n)}(\rdel^{\alpha_1},\ldots,\rdel^{\alpha_n}).$$
Let $\underline\alpha=(\alpha_0,\ldots,\alpha_n)$ and
$|\underline\alpha|=\sum_{i=0}^n \alpha_i.$ First, we will prove by induction on $q=|\underline\alpha|$
that
\begin{equation}\label{stablecomp}
\gamma(\underline\alpha\circ_k \underline\beta)=
\gamma(\underline\alpha)\circ_k \gamma(\underline\beta)
\end{equation}
If q=0 and  $\beta_0=0$ then
\begin{align*}
\gamma(\underline\alpha\circ_k\underline\beta)&=
\gamma(0,\ldots,0,\beta_1,\ldots,\beta_m,0,\ldots,0)\\
&=\mu(1,\ldots,1,\rdel^{\beta_1},\ldots,\rdel^{\beta_m},1,\ldots,1)\\
&=\mu\circ_k\mu(\rdel^{\beta_1},\ldots,\rdel^{\beta_m})\\
&=\gamma(\underline\alpha)\circ_k\gamma(\underline\beta)\\
\end{align*}
If $q=0$ and $\beta_0=1$, we use the notation of Proposition
\ref{homcyclicsphere},
where $\phi$ denotes $\sum_{i=1}^m \beta_i e_{i+k-1}$.  Thus,
\begin{align*}
\gamma(\underline\alpha\circ_k\underline\beta)&=
\sum_{l=0}^{k-1} \gamma(e_l+\phi)+ \sum_{l=k+1}^n
(-1)^{|\underline\beta|}\gamma(e_{l+m-1}+\phi)\\
&=(\ldel\circ_1\mu+\sum_{l=1,l\not=k}^{n}
\mu \circ_l\rdel)\circ_k\mu(\rdel^{\beta_1},\ldots,\rdel^{\beta_m})\\
&\stackrel{(\ref{equaldel})}=(\mu\circ_k\ldel)\circ_k
\mu(\rdel^{\beta_1},\ldots,
\rdel^{\beta_m})\\
&=\mu\circ_k(\ldel
\mu(\rdel^{\beta_1},\ldots,
\rdel^{\beta_m}))\\
&=\gamma(\underline\alpha)\circ_k\gamma(\underline\beta)
\end{align*}
We inductively assume (\ref{stablecomp}) is true for $q<p$ and
prove it for
$q=p$.

\begin{description}
\item[Case 1] If $q=1$, with $\alpha_0=1$, then
\begin{align*}
\gamma(\underline\alpha\circ_k\underline\beta)&=
\gamma(((1,0)\circ_1 (0,\ldots,0))\circ_k \underline\beta)\\
&=\gamma((1,0)\circ_1 ((0,\ldots,0)\circ_k \underline\beta)).
\end{align*}
It is therefore sufficient to prove (\ref{stablecomp}) for
$\underline\alpha=(1,0)$
and any $\underline\beta$. If $\beta_0=0$ it is trivial and
if $\beta_0=1$, then the putative relation (\ref{stablecomp}) vanishes
by the relation $\ldel\circ_1\ldel=0$.

\item[Case 2] Assume that $q\geq 1$ and that  if
$\alpha_0=1$ then $q\geq 2$. Thus, there exists $l\not=0$
such that $\alpha_l=1$. Since the cases $l<k$ and $l>k$ are symmetric,
we can assume first that $l<k$. In this case,
\begin{align*}
\gamma(\underline\alpha\circ_k\underline\beta)&=
\gamma(((\alpha_0,\ldots,\alpha_{l-1},0,\alpha_{l+1},\ldots,\alpha_n)
\circ_l(0,1))\circ_k\underline\beta) \\
&=(-1)^{|\underline\beta|}
\gamma(((\alpha_0,\ldots,\alpha_{l-1},0,\alpha_{l+1},\ldots,\alpha_n)
\circ_k\underline\beta)\circ_l (0,1))
\end{align*}
which proves, combined with the induction hypothesis, that it is
equivalent to prove (\ref{stablecomp}) for any $\underline\alpha$ and
$\underline\beta=(0,1)$. A straight-forward computation gives the
result
(using the identity $\rdel\circ_1\rdel=0$, in case $\alpha_k$=1).

In case $l=k$, we have by the induction hypothesis
\begin{align*}
\gamma(\underline\alpha\circ_k\underline\beta)&=
\gamma((\alpha_0,\ldots,\alpha_{k-1},0,\ldots,\alpha_{k+1},\ldots,\alpha_n)
\circ_k ((0,1)\circ_1 \underline\beta))\\
&=\gamma((\alpha_0,\ldots,\alpha_{k-1},0,\ldots,\alpha_{k+1},\ldots,\alpha_n)
\circ_k\gamma((0,1)\circ_1 \underline\beta)) \\
&=\gamma(\underline{\tilde\alpha})\circ_k
\gamma((0,1)\circ_1\underline\beta).
\end{align*}
On the other hand,
\begin{align*}
\gamma((0,1)\circ_1\underline\beta)&=
\sum_{l=1}^m \pm (-1)^{\beta_0}(\beta_l+1)[2]
\gamma(\beta_0,\beta_1,\ldots,\beta_l+1,\ldots,\beta_m)\\
&=\sum_{l=1}^{m}\pm (-1)^{\beta_0}
\ldel^{\beta_0}\mu(\rdel^{\beta_1},\ldots,\rdel^{\beta_l+1},\ldots,
\rdel^{\beta_m})\\
&=\sum_{l=1}^{m}(-1)^{\beta_0}\ldel^{\beta_0}
(\mu\circ_l\rdel)(\rdel^{\beta_1},\ldots,\rdel^{\beta_m})\\
&= (-1)^{\beta_0}\ldel^{\beta_0}\rdel
\mu(\rdel^{\beta_1},\ldots,\rdel^{\beta_m})\\
&= \ \
\rdel\ldel^{\beta_0}\mu(\rdel^{\beta_1},\ldots,\rdel^{\beta_m}) \\
&=\gamma(0,1)\circ_1 \gamma(\underline\beta),
\end{align*}
\noindent where the second-to-last equality follows from \ref{equardel},
and the next-to-last from the identity $\ldel\rdel=-\rdel\ldel$.

As a consequence
\begin{align*}
\gamma(\underline\alpha\circ_k\underline\beta)&=
(\gamma(\underline{\tilde\alpha})\circ_k
\gamma(0,1))\circ_k\gamma(\underline\beta)\\
&=\gamma(\underline\alpha)\circ_k\gamma(\underline\beta),
\end{align*}
\end{description}

The last step is to prove that $\gamma\circ\psi=Id$ and
$\psi\circ\gamma=Id$. The first equality holds for
the generators $\ldel, \rdel$ and $\mu$ and hence holds in general. The
second equality is proved again by induction on $|\underline\alpha|$~: if
$|\underline\alpha|=0$, then
$\psi(\gamma(\underline\alpha))=\psi(\mu)=\underline\alpha$.
Notice that $\psi(\gamma(0,1))=(0,1)$ and
$\psi(\gamma(1,0))=(1,0)$. Since any $\underline\alpha$ such that
$|\underline\alpha|>0$ is a composition $(1,0)\circ_1\underline\beta$
or
$\underline\beta\circ_i (0,1)$ with
$|\underline\beta|<|\underline\alpha|$
and since both $\gamma$ and $\psi$ are operad morphisms
(see (\ref{stablecomp})),
the induction hypothesis combined with the previous remark
yields the desired result.\endproof

\subsubth{Corollary} \label{corcyclicsphere}\sl An algebra $A$ over the operad $H_*(\lr)$
is a differential graded commutative and associative algebra
$(A,\cdot,\rdel)$,together with a differential $\ldel$  of degree $1$,
which anticommutes with $\rdel$. These operators satisfy the relation

$$\ldel(a)\cdot b=\ldel(a\cdot b) + (-1)^{|a|} a\cdot\rdel (b),\
\forall a,b\in A.$$
\rm

\subsubth{Remark} From the relations, it follows that
$\ldel\rdel$ is also a de\-ri\-vation of degree two.

\subsection{The family $\rr_\lambda$ of non cyclic operads on
$(S^1)^{n+1}$}

\subsubth{Definition} Let $\rr_\lambda$ be the collection of
spaces \hfill\break
$\rr_\lambda(n):= (S^{1})^{ (n+1)}$ together
with the following compositions:\\
$\circ_i: \rr_\lambda(n) \times \rr_\lambda(m) \rightarrow
\rr_\lambda(n+m-1)$
\begin{multline*}
(\theta_0,\dots,\th_n)\circ_i(\th'_0, \ldots, \th'_m) :=\\
(\th_0, \ldots, \th_{i-1},
\th'_1+ \th_i+\lambda \th'_{0}, \ldots, \th'_m +\th_i+\lambda\th'_{0},
\th_{i+1},\dots, \th_n)
\end{multline*}

Arguing in analogy to the earlier proofs, we obtain:

\subsubth{Proposition}\sl  $\rr_\lambda$ is a topological operad.
An algebra over $H_*(\rr_\lambda)$ is a differential
graded commutative and associative algebra $(A,\cdot,\rdel)$ together
with
an operator $\Delta : A\rightarrow A$ of degree $0$ and an operator
$\ldel : A\rightarrow A$ of degree 1 satisfying the relations
\begin{eqnarray*}
\Delta^2=\Delta, & \rdel\Delta=\Delta\rdel=\rdel,& \ldel\Delta=\ldel ,
\\
\Delta(ab)=&\Delta(a)b=a\Delta(b)&=ab ,\\
&\Delta\ldel=\lambda \rdel.&
\end{eqnarray*}
\rm

\subth{Remark} We have a morphism of operads $\rr_0 \rightarrow d$
which associates to $(\alpha_0,\alpha_1,\ldots,\alpha_n)$ the element
$(\alpha_1,\ldots,\alpha_n)$. Hence a differential graded commutative and
associative algebra $(A,\cdot,\del)$ is also a $H_*(\rr_0)$ algebra by
setting $\Delta=Id$, $\ldel=0$ and $\rdel=\del$.

Furthermore the map 
\begin{eqnarray*}
(S^1)^n &\rightarrow& (S^1)^{n+1} \\
(\theta_1, \dots, \theta_n) &\mapsto& (0,\theta_1, \dots, \theta_n) 
\end{eqnarray*}
yields an embedding of operads of $d \rightarrow \lr$
and  $d \rightarrow \rr_{\l}$.

\subth{Remark} 
The operads studied here have natural
extensions to more general situations.
For the operads $d$ and $\lr$, one can replace $S^1$ with any monoid $S$.
This fact lends itself to define a cyclic 
semi--direct product of a cyclic operad with a monoid in the spirit of
\cite{SW,MSS}.
The operad $\Cal Q$ is a special version of the operad of spaces \cite{K}.
To define the analogue of the operad $\rr_{\l}$ in this setting, one requires a monoid that is
a module over another monoid.

\end{document}